\documentclass[11pt,leq no,]{amsart}
\pdfoutput=1
\usepackage[utf8]{inputenc}
\usepackage[english]{babel}
\usepackage[T1]{fontenc}
\usepackage{fourier}
\usepackage{microtype}
\usepackage[,colorlinks=true,citecolor= DarkBlue,linkcolor=black]{hyperref}
\usepackage{times,framed,enumerate,stackrel}
\usepackage[svgnames,table]{xcolor}
\newcommand{\widebar}[1]{\mkern 1.5mu\overline{\mkern-1.5mu#1\mkern-1.5mu}\mkern 1.5mu}
\def\mathbi#1{\textbf{\textit #1}}
\usepackage{amssymb,amsbsy,amsfonts,amssymb,amscd,mathrsfs,amsmath,xspace,bm}
\def\[#1\]{\begin{align*}#1\end{align*}}
\newcommand{\pour}[1]{\qquad{\scriptstyle[#1]}}

\numberwithin{equation}{subsection}

\newcommand{\C}{\mathbb{C}\xspace}  \renewcommand{\P}{\mathbb{P}\xspace} \newcommand{\N}{\mathbb{N}\xspace} \newcommand{\Z}{\mathbb{Z}\xspace} \newcommand{\Q}{\mathbb{Q}\xspace}
\renewcommand{\O}{\mathcal{O}\xspace}
\renewcommand{\leq}{\leqslant} \renewcommand{\geq}{\geqslant}
\newcommand{\bydef}{\mathrel{\mathop:}=} \newcommand{\defby}{=\mathrel{\mathop:}}
\DeclareMathOperator{\moinsun}{-1}
\newcommand{\smbullet}{{\scriptscriptstyle \bullet}}
\newcommand{\abs}[1]{\lvert#1\rvert}
\usepackage{ifmtarg}
\makeatletter
\newtheoremstyle{myplain}
{12pt}{12pt}{\itshape}{}{\bfseries}{.}{ }
{\thmnumber{(\textup{#2})} \thmname{#1}\@ifmtarg{#3}{}{\,(\mdseries\thmnote{#3})}}
\newtheoremstyle{notitle}
{12pt}{12pt}{\itshape}{}{\bfseries}{}{ }
{\thmnumber{(\textup{#2})} \thmname{#1}\@ifmtarg{#3}{}{\,(\mdseries\thmnote{#3})}}
\makeatother
\theoremstyle{myplain}
\newtheorem*{THM}{Main Theorem}
\theoremstyle{notitle}
\newtheorem{Statement}[equation]{}
\author[Darondeau]
{Lionel~Darondeau}
\title
[Algebraic degeneracy of entire curves in \ensuremath{\P^{n}\setminus H}]
{
  Effective algebraic degeneracy of entire curves\\ 
  in complements of smooth projective hypersurfaces.\\
  {\tiny(Preliminary version)}
}
\date{6 Feb 2014}

\begin{document}
\maketitle
\section*{Introduction} 
In this work, the effective algebraic degeneracy of non-constant entire holomorphic curves having values in the complement of a generic smooth projective hypersurface of sufficiently high degree is established.
Namely we prove the:
\begin{THM}
  [Algebraic degeneracy of entire curves]
  \label{thm:main}
  If \(H\subset\P^{n}\) is a generic smooth projective hypersurface (of dimension \(n-1\)),
  having degree:
  \[
    d
    \,\geq\,
    (5n)^{2}\,n^{n},
  \]
  then there exists a proper subvariety \(Z\subset \P^{n}\), of codimension at least two, such that the image of every non constant entire curve
  \(
  f\colon \C \to\bigl(\P^{n}\setminus H\bigr)
  \)
  having values in the complement of \/\(H\), actually lies in \(\bigl(Z\setminus H\bigr)\).
\end{THM}

\medskip

In dimensions \(n=2\) \cite{MR1988702,MR2552951} and \(n=3\) \cite{MR2383820}, this theorem is already known with more precise lower bounds.
In \cite{MR2552951,MR1988123}, Erwan Rousseau also treat the case of hypersurfaces with several smooth components, in dimension \(2\).
In arbitrary dimension, the case of projective hypersurface with several smooth components will be treated in a forthcoming separate work \cite{D3}. 
Thanks to a lifting argument the problem reduces to the algebraic degeneracy of entire curves in complete intersections in \(P^{N}\), that can be treated with the same techniques that appear in this paper. These complete intersections were already studied in \cite{COM:9073852}.

\medskip

In order to prove the main theorem, we use the general strategy nicely described in \cite{MR2593279,div-rou-survey,paun-survey}, already successfully implemented in \cite{Paun2008,MR2383820,MR2331545,MR2593279}, that combines the \textsl{extrinsic approach} of Siu \cite{MR875091,MR958594,MR1420353,MR2077584,MR2543663,merkernew} with the \textsl{intrinsic approach} of Demailly \cite{MR1492539,MR1824906,MR2495771,MR1988702}.

The paper is organized as follows.
In the first section, the main result on degeneracy of entire curves drawn on complements of generic smooth projective hypersurfaces is proven.
The main idea is to produce \emph{many} algebraic differential equations satisfied by every non constant entire curve. 
Admitting for the moment the existence of \emph{one} such differential equation, and working in the universal family of complements of projective hypersurface of degree \(d\), lots of \emph{algebraically independent} differential equations are produced by using \textsl{low pole order slanted vector fields}~(\cite{MR2077584,MR2543663}).
These vector fields will be studied in our logarithmic context in the separate work \cite{DLPO}, which will generalize the dimension \(2\) and dimension \(3\) cases due to Erwan Rousseau (\cite{MR2383820,MR2552951}). 

The second section is devoted to proving the existence of one differential equation (result admitted in the first section). 
The intrinsic strategy outlined in~\cite{MR2495771} is followed. It relies on the construction of \textsl{invariant jet differentials} of Demailly~(\cite{MR1492539}), generalized to the logarithmic setting by Dethloff and Lu~(\cite{MR1824906}).
The use of weak algebraic Morse inequalities (\cite{MR1339712,MR1492539}) provides a control of the cohomology, sufficient for our goal.
It reduces the problem to the positivity of a certain intersection product \(I(d)\) on the \(k\)-th level of the \textsl{Demailly tower} of projectivized jet bundles.
This intersection product depends on parameters \(\underline{a}\) and \(\delta\), that can be adjusted later.
A variant of the multivariate residue formula of Berczi~(\cite{arxiv:1011.4710}), presented in \cite{D2013}, allows to integrate along the fibers of the Demailly tower. At this point, it remains to estimate one coefficient in the complicated Cauchy product of many (convergent) multivariate formal series.

In the third section, the computation of this coefficient is implemented. 
The overall approach, already adopted by Berczi in \cite{arxiv:1011.4710}, is to identify some central terms among the numerous combinations contributing to the coefficient of each power \(d^{p}\).
This identification leads to a modified version \(\widetilde{I}(d)\) of the polynomial \(I(d)\) having much more simple coefficients. It is easy to compute its largest root \(\widetilde{\lambda}\). 
Then, under suitable numerical hypotheses on the parameters \(\underline{a}\) and \(\delta\) ---\,that yield the effective degree bound \(d\geq(5n)^{2}n^{n}\) of the main theorem\,---, the estimation of the largest root of the polynomial \(I(d)\) is brought back to the much more easy computation of the largest root \(\widetilde{\lambda}\) of the simplified polynomial \(\widetilde{I}(d)\).
Many technical results are postponed to the appendices.

The last two sections: appendices \ref{apx:leadingCoeff} and \ref{apx:otherCoeffs}, form the \textsl{technical core} of this work.
In the first appendix, the leading coefficient of the intersection product \(I(d)\) is studied, and its positivity for a suitable (explicit) choice of the parameters \(\underline{a}\) and \(\delta\) is stated.
In the second appendix, the remaining coefficients of \(I(d)\) are studied and an upper bound \(\lambda\) for the largest root of this polynomial is derived.

\subsection*{Acknowledgments}
I am very thankful to Simone Diverio and Erwan Rousseau for discussing about \cite{MR2593279} and for regularly enquiring about the progress of this project. I am specially thankful to my thesis advisor, Joël Merker.

\section{Effective algebraic degeneracy of entire curves}
The formalism of jets, that will be now recalled, allows a coordinate-free description of differential equations and as a result more compact statements, although coordinates stay an essential tool in the machinery of proofs.

\subsection{Logarithmic jet bundles}
In what follows, \(\Delta\) denotes a complex disk of any radius, that can vary. Let \(\widebar{X}\) be a complex manifold.
\subsubsection*{Jet manifold}
Classically, the \textsl{\(k\)-jet manifold} of \(\widebar{X}\) is a coordinate-free construction with the same information as the \(k\)-th-order Taylor polynomial of germs of holomorphic maps \(\Delta\to\widebar{X}\).
For any germ of holomorphic map \(f\colon\Delta\to\widebar{X}\), the \textsl{\(k\)-jet} of \(f\) is:
\[
  f_{[k]}
  \bydef
  \left\{
    g\colon\Delta\to\widebar{X}
    \;\middle\vert\;
    \text{\(g\) and \(f\) osculate to order \(k\) at the origin}
  \right\}.
\]
Then, for a target point \(x\in\widebar{X}\), introduce:
\[
  J_{k}\widebar{X}_{x}
  \bydef
  \left\{
    f_{[k]}
    \;\middle\vert\;
    f\colon(\Delta,0)\to(\widebar{X},x)
  \right\}.
\]
The collection of these spaces is a fiber bundle, called the \textsl{\(k\)-jet manifold} of \(\widebar{X}\):
\[
  J_{k}\widebar{X}
  \bydef
  \bigcup_{x\in\widebar{X}}
  J_{k}\widebar{X}_{x}
  \;\stackrel{\eta}{\to}\;
  \widebar{X},
\]
where the map \(\eta\) is the evaluation of the jets at the origin:
\[
  \eta\bigl(f_{[k]}\bigr)=f(0).
\]
A jet \(f_{[k]}\in J_{k}\widebar{X}\) is termed \textsl{singular}, if \(f\) is stationary (\textit{i.e.} if \(f'(0)=0\)), and \textsl{regular} otherwise.

For a meromorphic differential 1-form \(\omega\in\Gamma_{\text{loc}}\bigl(T_{\widebar{X}}^{\star}\bigr)\), and a germ of holomorphic map \(f\colon\Delta\to\widebar{X}\), the pullback \(f^{\star}\omega\) is necessarily of the form:
\[
  f^{\star}\omega
  =
  A(t)dt,
\]
for a meromorphic function \(A\colon\C\to\C\).
Thus, each such \(1\)-form \(\omega\) induces a meromorphic map:
\[
  \widetilde{\omega}
  \colon
  f_{[k]}
  \mapsto
  \bigl(
  A(t),
  A'(t),
  \dotsc,
  A^{(k-1)}(t)
  \bigr)
  \pour{
    J_{k}X\to\C^{k}
  }.
\]
On an open set \(U\subset\widebar{X}\), the trivialization associated to a meromorphic local coframe \(\omega_{1}\wedge\dotsb\wedge\omega_{n}\neq0\) is:
\[
  \left\{%}
  \begin{array}{ccc}
    \Gamma_{U}\bigl(J_{k}\widebar{X}\bigr)
    &\to 
    &U\times\bigl(\C^{k}\bigr)^{n}
    \\
    \sigma
    &\mapsto 
    &
    \eta\circ \sigma
    \,;\,
    \widetilde{\omega}_{1}\circ \sigma,
    \dotsc,
    \widetilde{\omega}_{n}\circ \sigma
  \end{array}
  \right..
\]
The components \(A_{i}^{(j)}\) of \(\widetilde{\omega}_{i}\circ\sigma\) are called \textsl{jet coordinates} and correspond to the derivatives of the germs \(\sigma_{x}\colon t\mapsto \sigma_{x}(t)\) with respect to the complex variable \(t\in\Delta\).

\subsubsection*{Logarithmic jet manifold along a normal crossing divisor}
A divisor \(D\subset\widebar{X}\) has only \textsl{normal crossings} if at each point \(x\in\widebar{X}\), there is an integer \(l=l(x)\) and a centered coordinate system \(z_{1},\dotsc,z_{l},z_{l+1},\dotsc,z_{n}\) on \(\widebar{X}\) around \(x\) such that:
\[
  D
  \simeq
  \mathrm{div}(z_{1}\dotsm z_{l})
  \subset
  \C^{n}.
\]
For such a normal crossing divisor \(D\), one defines:
\[
  \mathcal{T}^{\star}_{\widebar{X}}(\log D)_{x}
  \bydef
  \O_{\widebar{X}}\frac{dz_{1}}{z_{1}}
  +\dotsb+
  \O_{\widebar{X}}\frac{dz_{l}}{z_{l}}
  +
  \O_{\widebar{X}}dz_{l+1}
  +\dotsb+
  \O_{\widebar{X}}dz_{n}.
\]
It is a locally free \(\O_{\widebar{X}}\)-module of rank \(n\) that is the sheaf of sections of a vector bundle \(T^{\star}_{\widebar{X}}(\log D)\), called the \textsl{logarithmic cotangent bundle of \(\widebar{X}\) along \(D\)} (\cite{MR637060}).
\medskip

A local section \(\sigma\in\Gamma_{U}(J_{k}\widebar{X})\) over an open set \(U\in\widebar{X}\) is termed \textsl{logarithmic} if for any logarithmic cotangent vector field \(\omega\in\Gamma_{V}(T^{\star}_{\widebar{X}}(\log D))\), defined on a smaller open set \(V\subset U\), the obtained meromorphic function \(\widetilde{\omega}\circ\sigma\) is actually holomorphic, or, in other words, if \(\sigma\) has holomorphic jet coordinates in the adapted logarithmic coframe generating \(\mathcal{T}^{\star}_{\widebar{X}}(\log D)\). These sections define a vector subbundle:
\[
  J_{k}\widebar{X}(-\log D)
  \subset
  J_{k}\widebar{X},
\]
called the \textsl{logarithmic \(k\)-jet manifold of \(\widebar{X}\) along \(D\)} (\cite{MR859200}).

The subset of \textsl{singular logarithmic \(k\)-jets} is defined as the Zariski closure of the set of singular \(k\)-jets in the logarithmic \(k\)-jet manifold:
\[
  J_{k}^{\text{sing}}\widebar{X}(-\log D)
  \bydef
  \overline{
    J_{k}^{\text{sing}}\widebar{X}
    \cap
    J_{k}\widebar{X}(-\log D).
  }^{\text{Zar}}.
\]
The logarithmic \(k\)-jets that are not singular are said \textsl{regular}.

\subsubsection*{Jet differentials}
Recall the concept of jet differentials, after \cite{MR609557,MR1492539,MR1824906}. ``It is a coordinate-free description of the holomorphic differential equations that a germ of curve may satisfy''.

One has a \(\C^{\star}\)-action on \(J_{k}\widebar{X}(-\log D)\) by rescaling of the source.
Indeed if \(h_{\lambda}\) is the homothety with ratio \(\lambda\in\C^{\star}\), and \(f\) is a logarithmic jet along \(D=\mathrm{div}(z_{1}\dotsm z_{l})\), the jet coordinates change as follows:
\[
  \begin{cases}
    \bigl(\log\circ f_{i}\circ h_{\lambda}\bigr)^{(j)}
    =
    \lambda^{j}\;\bigl(\log\circ f_{i}\bigr)^{(j)}\circ h_{\lambda}
    &\pour{i=1,\dotsc,l}
    \\
    \bigl(f_{i}\circ h_{\lambda}\bigr)^{(j)}
    =
    \lambda^{j}\;f_{i}^{(j)}\circ h_{\lambda}
    &\pour{i=l+1,\dotsc,n}
  \end{cases}.
\]

The Faà di Bruno formulae show that the concept of polynomial on the fibers of \(J_{k}\widebar{X}(-\log D)\) makes sense.
One can thus consider the \textsl{Green-Griffiths jet bundle} of differential operators of \textsl{order} \(k\) and \textsl{weighted degree} \(m\):
\[
  \mathcal{E}_{k,m}T^{\star}_{\widebar{X}}(\log D\bigr)
  \to
  \widebar{X},
\]
the fiber of which consists of the complex valued polynomials \(Q(f',\dotsc,f^{(k)})\) on the fibers of \(J_{k}\widebar{X}(-\log D)\) of weighted degree \(m\) with respect to the \(\C^{\star}\)-action by rescaling of the source, that is:
\[
  Q\bigl(
  (f\circ h_{\lambda})_{[k]}
  \bigr)
  =
  \lambda^{m}\,Q\bigl(f_{[k]}\bigr).
\]
A nice way to describe the polynomials \(Q\) is to say that in any monomials of \(Q\), \(m\) is the ``number of primes'' appearing .

\subsection{Fundamental vanishing theorem}
The link between algebraic degeneracy and jet bundles is nowadays classical. On has indeed the following \textsl{fundamental vanishing theorem} (\cite{MR1492539,MR2364143})
\begin{Statement}
    %[Ahlfors-Schwarz lemma]
  \label{thm:ahlfors_schwartz}
  Let \(P\) be a non zero global jet differential of order \(k\) and weighted degree \(m\),
  vanishing with order \(\varepsilon>0\):
  \[
    P
    \in
    H^{0}
    \bigl(
    \P^{n},
    \mathcal{E}_{k,m}T_{\P^{n}}^{\star}(\log H)
    \otimes 
    \O_{\P^{n}}(\varepsilon)^{\vee}
    \bigr)
    \setminus
    \big\{0\big\},
  \] 
  then every non constant entire holomorphic curve \(f\colon\C\to\bigl(\P^{n}\setminus H\bigr)\) must satisfy
  the corresponding algebraic differential equation of degree \(k\) with \(m\) ``primes'':
  \[
    P\bigl(f^{(1)},\dotsc,f^{(k)}\bigr)
    \equiv
    0.
  \]
\end{Statement}
\subsection{Universal family of complements of projective hypersurfaces}
Let \( S\) be the space of parameters of degree \(d\) projective hypersurface of \(\P^{n}\):
\[
  S
  \bydef
  \P\, H^{0}\bigl(\P^{n},\O_{\P^{n}}(d)\bigr).
\]
Consider the universal family of degree \(d\) hypersurfaces:
\[
  \mathcal{H}
  \bydef
  \{(x,P)\in\P^{n}\times S\colon P(x)=0\},
\] 
and the complement \(\mathcal{X}\) of this universal family:
\[
  \mathcal{X}
  \bydef
  \bigl(
  \P^{n}\times S
  \bigr)
  \setminus
  \mathcal{H}.
\]
Denote the two natural projections to the factors of \(\P^{n}\times S\) by:
\[
  \P^{n}
  \stackrel{\mathrm{pr}_{1}}\longleftarrow
  \P^{n}\times S
  \stackrel{\mathrm{pr}_{2}}\longrightarrow
  S.
\]

If \(s\in S\) is a point of the parameter space of \(\mathcal{H}\), then:
\[
  \mathcal{X}_{s}
  \bydef
  (\mathrm{pr}_{2}\vert_{\mathcal{X}})^{\moinsun}(s)
  =
  \bigl(\P^{n}\setminus H_{s}\bigr)\times\{s\},
\]
where \(H_{s}\subset\P^{n}\)
is the projective hypersurface of degree \(d\) parametrized by \(s\):
\[
  H_{s}
  \bydef
  \mathrm{pr}_{1}
  \bigl(
  (\mathrm{pr}_{2}\vert_{\mathcal{H}})^{\moinsun}
  (s)
  \bigr).
\] 
Moreover, for a generic \(s\in S\), the hypersurface \(H_{s}\) is smooth.

\subsubsection*{Slanted vector fields}
The universal family \(\mathcal{H}\) is a %smooth(?) 
normal crossing divisor of degree \((d,1)\) in \(\P^{n}\times S\). 
The space of the logarithmic \(k\)-jets along \(\mathcal{H}\) of the log-manifold \(\bigl(\P^{n}\times S,\mathcal{H}\bigr)\) will be denoted, as above, by:
\[
  J_{k}(\P^{n}\times S)(-\log\mathcal{H}).
\]

The \textsl{vertical jet fields} are the fields tangent to the fibers of \(\mathrm{pr}_{2}\colon\P^{n}\times S\to S\).
In other words, let:
\[
  \mathrm{pr}_{2\star}^{[k]}
  \colon
  J_{k}(\P^{n}\times S)
  \to
  J_{k}( S)
\] 
be the extension to \(k\)-jets of the bundle morphism:
\[
  {\mathrm{pr}_{2}}_{\star}
  \colon
  T(\P^{n}\times S)
  \to
  T( S),
\] 
then:
\[
  J_{k}^{\mathrm{vert}}(\P^{n}\times S)
  \bydef
  \ker\;
  \mathrm{pr}_{2\star}^{[k]}
  \subset
  J_{k}(\P^{n}\times S).
\]

The submanifold of \textsl{vertical regular logarithmic jet fields} of order \(k\):     
\[
  J_{k}^{\mathrm{vert,reg}}(\P^{n}\times S)(-\log\mathcal{H})
  \bydef
  J_{k}^{\mathrm{reg}}(\P^{n}\times S)(-\log\mathcal{H})
  \cap
  J_{k}^{\mathrm{vert}}(\P^{n}\times S),
\]
consists of regular logarithmic jets tangent to the fibers of the natural projection \(\mathrm{pr}_{2}\) over the parameter space. 

  Recall that \(\eta\) denotes the natural projection 
  (\textit{i.e.} the evaluation of the jets):
  \[
    \eta
    \colon
    J_{k}^{\mathsf{vert}}(\P^{n}\times S)(-\log\mathcal{H})
    \to
    \bigl(\P^{n}\times S\bigr).
  \] 
In a separate work \cite{DLPO}, we prove the following statement:
\begin{Statement}
    %[Low pole order slanted vector frames]
  \label{thm:lpo_vector_frames}
  For degrees \(d>k\), the (twisted) holomorphic tangent sheaf to vertical jets of the log manifold \((\P^{n}\times S,\mathcal{H})\):
  \[
    T_{J_{k}^{\mathsf{vert}}(\P^{n}\times S)(-\log\mathcal{H})}
    \otimes
    \eta^{\star}
    \left(
    \O_{\P^{n}}\bigl(k(k+2)\bigr)
    \otimes
    \O_{ S}(1)
    \right)
  \]
  is generated by its global holomorphic sections over the subspace:
  \[
    J_{k}^{\mathsf{vert,reg}}(\P^{n}\times S)(-\log\mathcal{H})
    \setminus 
    \eta^{\moinsun}\mathcal{H},
  \]
  of regular vertical logarithmic \(k\)-jets of holomorphic curves avoiding \(\mathcal{H}\). 
\end{Statement}

\subsection{Proof of the main result}
A strategy~\cite{MR2077584,Paun2008,MR2331545,MR2383820,MR2593279,div-rou-survey} for proving the \textsl{algebraic degeneracy} of entire curves, \textit{i.e.} for the obtainment of  an algebraic differential equation of degree \(0\) that every entire curve has to satisfy is then to generate many global jet differentials. 
Using the above vanishing theorem~\eqref{thm:ahlfors_schwartz}, this gives algebraic differential equations that every entire curve must satisfy.
Next, one can (morally) get rid of the differentials \(f^{(1)},\dotsc,f^{(k)}\) by algebraic elimination.

\subsubsection*{Existence of algebraic differential equations}

In \S\ref{sec:existence}, using algebraic Morse inequalities, we prove the following statement:
\begin{Statement}
  \label{thm:existence}
  For any smooth hypersurface \(H\subset\P^{n}\) of degree \(d\) and any positive rational number \(\delta\) such that:
  \[
    d\geq52n^{n}
    \text{ and }
    35n^{n}\delta\leq1,
  \]
  the vector space of logarithmic global jet differentials along \(H\) of order \(n\) and weighted degree \(m\gg d\) having values in the dual of the ample line bundle \(\mathcal{K}_{\P^{n}}(H)^{2\delta m}\) has positive dimension:
  \[
    \dim\,
    H^{0}
    \Bigl(
    \P^{n},
    \mathcal{E}_{n,m}T_{\P^{n}}^{\star}\bigl(\log H\bigr)
    \otimes 
    \bigl(\mathcal{K}_{\P^{n}}(H)^{\vee}\bigr)^{2\delta m}
    \Bigr)
    \geq
    1.
  \]
\end{Statement}

The sections of 
\(
\mathcal{E}_{n,m}T_{\P^{n}}^{\star}(\log H)
\otimes 
\mathcal{K}_{\P^{n}}(H)^{-2\delta m}
\)
can be interpreted as invariant maps \(J_{n}\P^{n}(-\log H)\to\eta^{\star}\mathcal{K}_{\P^{n}}(H)^{-2\delta m}\).
Recall moreover that:
\[
  \mathcal{K}_{\P^{n}}(H)
  =
  \O_{\P^{n}}(d-n-1).
\]

\medskip

Let \(f\colon\C\to\mathcal{X}\) be a non constant entire curve having values in the complement of the universal family of projective hypersurfaces of degree \(d\). 
Assume \(f\) is tangent to the logarithmic relative tangent bundle of \(\mathrm{pr}_{2}\):
\[
  \mathcal{V}(-\log\mathcal{H})
  \bydef
  \ker\bigl({\mathrm{pr}_{2}}_{\star}\bigr).
\] 
This assumption exactly means that \(\mathrm{pr}_{2}\circ f\equiv s_{0}\) is constant. 
In other words, \(f\) maps \(\C\) to the slice:
\[
  \mathcal{X}_{s_{0}}
  \bydef
  \bigl(\P^{n}\setminus H_{s_{0}}\bigr)
  \times
  \big\{s_{0}\big\}
  \subset 
  \mathcal{X}.
\]
Generically, one can assume that \(H_{s_{0}}\) is smooth.

By the very definition of the relative tangent bundle, one has:
\[
  \mathcal{V}(-\log\mathcal{H})\vert_{\mathcal{X}_{s}}
  =
  T_{\mathcal{X}_{s}}(-\log\mathcal{H}_{s})
  =
  (\mathrm{pr}_{2})^{\star}\,T_{\P^{n}}(-\log H_{s}).
\]

By the above statement \eqref{thm:existence}, 
for a suitable choice of \(m\gg d\gg1\) and \(0<\delta\ll1\), independently of \(s_{0}\),
there exists a non zero jet differential having positive vanishing order:
\[
  \sigma_{s_{0}}
  \in
  H^{0}\Bigl(
  \P^{n}\times\{s_{0}\},
  \mathcal{E}_{n,m}
  \mathcal{V}^{\star}(\log\mathcal{H})
  \otimes
  \eta^{\star}
  \O_{\P^{n}}(2\delta m(d-n-1))^{\vee}
  \big\vert_{\mathcal{X}_{s_{0}}}
  \Bigr),
\]
with zero locus:
\[
  Z_{s_{0}}
  \bydef
  \{
    x\in\P^{n}
    \colon
    \sigma_{s_{0}}(x) = 0
  \}
  \subsetneq
  \P^{n}\times\{s_{0}\}.
\]

Because the vanishing order of the section \(\sigma_{s_{0}}\) is positive, one can apply the fundamental vanishing theorem and obtain:
\(
\sigma_{s_{0}}
\bigl(f_{[n]}\bigr)
\equiv
0.
\) 

In fact, it will be established that the image of \(f\) lies in the common zero set \(Z_{s_{0}}\) of all holomorphic coefficients of the polynomial \(\sigma_{s_{0}}\left(f^{(1)},\dotsc,f^{(n)}\right)\).
The above differential equation is thus a meaningless artefact. Indeed, the derivative of \(f\) are actually not involved in the vanishing of \(\sigma_{s_{0}}\left(f_{[n]}\right)\).

\subsubsection*{Lie derivative}
In order to use Lie derivative, extend first \(\sigma_{s_{0}}\) as a holomorphic family of non zero jet differentials:
\[
  \sigma
  \bydef
  \Bigl\{
    \sigma_{s}
    \in 
    H^{0}\bigl(
    \P^{n}\times\{s\},
    \mathcal{E}_{n,m}T_{\P^{n}}^{\star}(\log H_{s})
    \otimes
    \O_{\P^{n}}(2\delta m(d-n-1))^{\vee}
    \bigr)
  \Bigr\}.
\]
By a classical extension result (\cite[(3.2)]{MR0463157}), there exists a Zariski closed subset \(\Sigma\subset S\) such that if \(s_{0}\) lies in \(S\setminus\Sigma\), then \(\sigma_{s_{0}}\) can be extended as a holomorphic family \(\sigma\) to a Zariski dense open set containing \(s_{0}\).
Make the two generic assumptions that \(s_{0}\notin\Sigma\) and that \(H_{s_{0}}\subset\P^{n}\) is smooth and, in order to get rid of the term \(\eta^{\star}\O_{S}(1)\) in \eqref{thm:lpo_vector_frames}, remove a supplementary hyperplane of \(S\), if necessary.

Working locally on a neighbourhood \(U\) of one such \(s_{0}\), consider the section:
\[
  \sigma
  \colon
  J_{n}\mathcal{V}(-\log\mathcal{H})
  \vert_{\mathrm{pr}_{2}^{\moinsun}U}
  \to
  \eta_{1}^{\star}
  \O_{\P^{n}}\bigl(
  -2\delta m(d-n-1)
  \bigr)
  \vert_{\mathrm{pr}_{2}^{\moinsun}U}.
\]
and \(p\) slanted vector fields provided by \eqref{thm:lpo_vector_frames}:
\[
  \mathsf{v}_{1},\dotsc,\mathsf{v}_{p}
  \colon
  \mathrm{pr}_{2}^{\moinsun}U
  \to
  T_{J_{n}\mathcal{V}(-\log\mathcal{H})}
  \otimes
  \eta^{\star}
  \O_{\P^{n}}(n(n+2))
  \vert_{\mathrm{pr}_{2}^{\moinsun}U}.
\]
These slanted vector fields are useful in order to generate lots of new vector fields by Lie derivative of the given section \(\sigma\).
Indeed, by induction, the Lie derivative
\(
\bigl(\mathsf{v}_{p}\dotsm\mathsf{v}_{1}\cdot\sigma\bigr)_{s_{0}}
\)
is a non zero section of the vector space:
\[
  H^{0}\Bigl(
  \P^{n}\times\{s_{0}\},
  \mathcal{E}_{n,m}T_{\P^{n}}^{\star}(\log H_{s_{0}})
  \otimes
  \O_{\P^{n}}\bigl(
  -2m\delta(d-n-1)+p\,n(n+2)
  \bigr)
  \Bigr).
\]
Each derivative decreases the vanishing order by at most \(n(n+2)\).

While the vanishing order of 
\(
\bigl(\mathsf{v}_{p}\dotsm\mathsf{v}_{1}\cdot\sigma\bigr)_{s_{0}}
\)
is still positive, 
\textit{i.e.} while:
\[
  2\delta m(d-n-1)
  >
  pn(n+2),
\]
the above argument remains valid, and one infers that:
\[
  \bigl(\mathsf{v}_{p}\dotsm\mathsf{v}_{1}\cdot\sigma\bigr)_{s_{0}}
  \bigl(f_{[n]}\bigr)
  \equiv 0.
\]
In this way, \emph{lots of algebraic differential equations} can be obtained, as soon as the vanishing order of the initial section \(\sigma\) is large enough.

\subsubsection*{Sufficient lower bound on the degree}
The element of contradiction is the following.
If \(f(\C)\not\subset Z_{s_{0}}\), pick \(t_{0}\in\C\) such that \(f(t_{0})\not\in Z_{s_{0}}\).
One may assume that \(f'(t_{0})\neq0\), because \(f_{[n]}(\C)\not\subset J_{n}^{\mathrm{sing}}\mathcal{V}(-\log\mathcal{H})\).
Now, by the result of global generation ({\cite[p.12]{MR2593279}}):
\begin{Statement}
  If \(f(t_{0})\not\in Z_{s}\) and \(f'(t_{0})\neq0\), then there exist \(q\leq m\) invariant vector fields \(\mathsf{v}_{1},\dotsc,\mathsf{v}_{q}\) such that:
  \[
    \bigl(\mathsf{v}_{q}\dotsm\mathsf{v}_{1}\cdot\sigma\bigr)_{s_{0}}
    \bigl(
    f_{[n]}(t_{0})
    \bigr)
    \neq
    0.
  \]
\end{Statement}
Thus, if for any \(p\leq m\) and any invariant vector fields \(\mathsf{v}_{1},\dotsc,\mathsf{v}_{p}\), one has the vanishing:
\[
  0
  \equiv
  \bigl(\mathsf{v}_{p}\dotsm\mathsf{v}_{1}\cdot\sigma\bigr)_{s_{0}}
  \bigl(f_{[n]}\bigr),
\] 
the sought contradiction is obtained.
In particular, by the fundamental vanishing theorem, it is the case if the vanishing order is positive for all \(p\leq m\). In other words, if the degree \(d\) satisfies:
\[
  2\delta m(d-n-1)
  >
  mn(n+2),
\]
then, the entire curves with values in the complement of a generic hypersurface of degree \(d\) are algebraically degenerate.
After dividing by \(m\), the condition on the degree becomes:
\[
  d
  >
  (n+1)
  +
  \frac{n(n+2)}{2\delta}.
\]

Now, recall that by assumption
\(
35n^{n}\,\delta
\leq
1.
\)
After fixing:
\(
35n^{n}\,\delta
=
1,
\)
in order to minimize the lower bound, one finally obtains the stronger condition (for \(n\geq5\)):
\[
  d
  \geq
  25\,n^{n+2}
  >
  n^{n+2}
  \left(
  \frac{n+1}{n^{n+2}}
  +
  \frac{35(n+2)}{2n}
  \right).
\]
\begin{sffamily}
  \small
  One can drop the hypothesis \(n\geq5\), because for \(n\leq5\), the explicit computation of the largest positive root, using computer algebra system, shows that the result is granted for a lower bound smaller than \(25n^{n+2}\).
\end{sffamily}

Notice that this bound \((5n)^{2}n^{n}>52n^{n}\) is sufficiently large in order to prove the existence of one global jet differential.

\subsubsection*{Codimension of the degeneracy locus}
In \cite{MR2746466}, the authors show that the degeneracy locus of entire curves has no divisorial components. As a result, its codimension is at least two.

\section{Existence of global logarithmic jet differentials}
\label{sec:existence}
\setcounter{equation}{0}
In this section, the intrinsic strategy outlined in~\cite{MR2495771} is followed in order to prove the existence theorem:
\begin{Statement}
  For any smooth hypersurface \(H\subset\P^{n}\) of degree \(d\) and any positive rational number \(\delta\) such that:
  \[
    d
    \geq
    (52\,n^{n})
    \text{ and }
    (35\,n^{n})\,\delta
    \leq
    1,
  \]
  the vector space of logarithmic global jet differentials along \(H\) of order \(k=n\) and weighted degree \(m\gg n^{n}\) vanishing on \(\mathcal{K}_{\P^{n}}(H)^{2\delta m}\) has positive dimension:
  \[
    \dim\,
    H^{0}
    \Bigl(
    \P^{n},
    \mathcal{E}_{n,m}T_{\P^{n}}^{\star}(\log H)
    \otimes 
    \mathcal{K}_{\P^{n}}(H)^{-2\delta m}
    \Bigr)
    \geq
    1.
  \]
\end{Statement}

A general intrinsic method for obtaining global jet differentials relies on the construction of Demailly~(\cite{MR1492539}), generalized to the logarithmic setting by Dethloff and Lu~(\cite{MR1824906}), that will be now described.

\subsection{Demailly tower of logarithmic directed manifolds}
A \textsl{log-directed manifold} is a triple \(\bigl(\widebar{X},D,V\bigr)\), where \(\bigl(\widebar{X},D\bigr)\) is a log manifold and \(V\) is a holomorphic subbundle of the logarithmic tangent bundle \(T_{\widebar{X}}(-\log D)\).
Given such a log-directed manifold \((\widebar{X}_{0},D_{0},V_{0})\), one constructs on \(\widebar{X}_{0}\) the \textsl{Demailly tower of log directed manifolds}:
\[
  (\widebar{X}_{\kappa},D_{\kappa},V_{\kappa})
  \stackrel{\pi_{\kappa}}\longrightarrow
  (\widebar{X}_{\kappa-1},D_{\kappa-1},V_{\kappa-1})
  \stackrel{\pi_{\kappa-1}}\longrightarrow
  \dotso
  \stackrel{\pi_{1}}\longrightarrow
  (\widebar{X}_{0},D_{0},V_{0}),
\]
by induction on \(\kappa\geq0\). 
This construction has the same formal properties as in the so-called \textsl{compact case}, \textit{i.e.} where there is no divisor \(D_{0}\). Here, in the genuine logarithmic setting, \(V_{k}\) is a holomorphic subbundle of the logarithmic tangent bundle \(T_{\widebar{X}_{k}}(-\log D_{k})\). 

Recall quickly the inductive step \(\bigl(\widebar{X}',D',V'\bigr)\stackrel{\pi}\to\bigl(\widebar{X},D,V\bigr)\) of the construction of the log Demailly tower (\cite{MR1824906}).
The space \(\widebar{X}'\) is the total space \(P(V)\) of the projective bundle of lines of \(V\):
\[
  \widebar{X}'
  \bydef
  P(V)
  \stackrel{\pi}\longrightarrow
  \widebar{X}.
\]
In order to make \(\pi\) a log-morphism it is natural to set:
\[
  D'
  \bydef
  \pi^{\moinsun}(D)
  \subset
  \widebar{X}'.
\]

The \textsl{relative tangent bundle} \(T_{\pi}\) of the log-morphism \(\pi\colon(\widebar{X},D)\to(\widebar{X}',D')\) is defined as the kernel of the differential \(\pi_{\star}\):
\begin{equation}
  \label{eq:relative_tangent_bundle}
  0
  \to
  T_{\pi}
  \hookrightarrow
  T_{\widebar{X}'}(-\log D')
  \stackrel{\pi_{\star}}
  \longrightarrow
  \pi^{\star}
  T_{\widebar{X}}(-\log D)
  \to
  0.
\end{equation}
Keep in mind that \(V\) is a subbundle of the logarithmic tangent bundle \(T_{\widebar{X}}(-\log D)\), displayed in the right of this short exact sequence, and that \(V'\) has to be a subbundle of the logarithmic tangent bundle \(T_{\widebar{X}'}(-\log D')\), displayed in the center of this short exact sequence.

The tautological line bundle of \(\widebar{X}'=P(V)\) is a subbundle of \(\pi^{\star}T_{\widebar{X}}(-\log D)\) because:
\[
  \O_{\widebar{X}'}(-1)
  \subset
  \pi^{\star}V
  \subset
  \pi^{\star}T_{\widebar{X}}(-\log D).
\]
One can thus define a subbundle \(V'\subset T_{X'}(-\log D')\), locally isomorphic to \(\pi^{\star}V\), by taking:
\[
  V'
  \bydef
  (\pi_{\star})^{\moinsun}
  \O_{\widebar{X}'}(-1).
\]
Equivalently, \(V'\) is defined by the following short exact sequence:
\begin{equation}
  \label{eq:v_prime}
  0
  \to
  T_{\pi}
  \hookrightarrow
  V'
  \stackrel{\pi^{\star}}\longrightarrow
  \O_{\widebar{X}'}(-1)
  \to
  0.
\end{equation}
Comparing the two short exact sequences \eqref{eq:relative_tangent_bundle} and \eqref{eq:v_prime}, notice that in the left, one keeps all the vertical directions whereas in the right, one keeps only the single ``tautological'' direction among all horizontal directions.

\subsection{The direct image formula}
For any integer \(i=1,\dotsc,\kappa-1\), the composition of the projections \(\pi_j\colon X_j\to X_{j-1}\) yields a projection:
\[
  \pi_{\kappa,i}\bydef\pi_{i+1}\circ\dotsb\circ\pi_{\kappa}
  \colon
  X_{\kappa}\to X_i.
\] 
One can pullback the tautological line bundle of \(X_i\) to the \(\kappa\)-th level. 
We obtain \(\kappa\) line bundles:
\[
  (\pi_{\kappa,1})^{\star}
  \O_{\widebar{X}_{1}}(-1),\;
  (\pi_{\kappa,2})^{\star}
  \O_{\widebar{X}_{2}}(-1),\;
  \dotsc\;,\;
  (\pi_{\kappa,\kappa-1})^{\star}
  \O_{\widebar{X}_{\kappa-1}}(-1),\;
  \O_{\widebar{X}_{\kappa}}(-1).
\] 
The linear combinations of these line bundles with integer coefficients \((a_{1},\dotsc,a_{\kappa})\in\Z^{\kappa}\):
\[
  \O_{\widebar{X}_{\kappa}}(a_{1},\dotsc,a_{\kappa})
  \bydef
  \O_{\widebar{X}_{\kappa}}(a_{\kappa})
  \otimes
  (\pi_{\kappa,\kappa-1})^{\star}\O_{\widebar{X}_{\kappa-1}}(a_{\kappa-1})
  \otimes
  \dotsb
  \otimes
  (\pi_{\kappa,1})^{\star}
  \O_{\widebar{X}_{1}}(a_{1}).
\]
enjoy positivity properties~(\cite{MR1492539}): 

In \cite{D2013}, we show that it is natural to work with the \(\kappa\) line bundles:
\[
  \mathcal{L}_{i}^{\vee}
  \bydef
  \O_{\widebar{X}_{i}}(1,\dotsc,1)
  \pour{i=1,\dotsc,\kappa}.
\]
For integers coefficients \((a_{1},\dotsc,a_{\kappa})\in\Z^{\kappa}\), in analogy with the notation \(\O_{\widebar{X}_{k}}(\underline{a})\), let us use the notation:
\[
  \mathcal{L}(a_{1},\dotsc,a_{\kappa})
  \bydef
  (\pi_{\kappa,1})^{\star}
  (\mathcal{L}_{1}^{\vee})^{a_{1}}
  \otimes
  \dotsb
  \otimes
  (\pi_{\kappa,\kappa})^{\star}
  (\mathcal{L}_{\kappa}^{\vee})^{a_{\kappa}}.
\]

Obviously, the translation formula is:
\begin{equation}
  \label{eq:translation_line_bundles}
  \mathcal{L}(a_{1},a_{2},\dotsc,a_{\kappa})
  =
  \O_{\widebar{X}_{\kappa}}(a_{1}+a_{2}+\dotsb+a_{\kappa},a_{2}+\dotsb+a_{\kappa},\dotsc,a_{\kappa}).
\end{equation}

Considering this formula \eqref{eq:translation_line_bundles} ---\,that can be thought of as a plain change of variables\,---,
one can easily transpose existing results on the line bundles \(\O_{\widebar{X}_{\kappa}}(\underline{a})\) (\cite{MR1492539,MR1824906})
into results on the line bundles \(\mathcal{L}(\underline{a})\).
One has the following \textsl{direct image formula}:
\begin{Statement} 
  For any \(\kappa\)-tuple \((a_{1},\dotsc,a_{\kappa})\in \Z^{\kappa}\) with:
  \[
    a_{i}+\dotsb+a_{\kappa}
    \geq
    0
    \pour{i=1,\dotsc,\kappa},
  \]
  the direct image
  \(
  (\pi_{\kappa,0})_{\star}
  \mathcal{L}(a_{1},\dotsc,a_{\kappa})
  \) 
  may be seen as a subsheaf of the Green-Griffiths bundle \(\mathcal{E}_{\kappa,\mu}V_{0}^{\star}(\log D_{0})\),
  where:
  \[
    \mu(\underline{a})
    \bydef
    1\,a_{1}+2\,a_{2}\dotsb+\kappa\,a_{\kappa}
    \in\N.
  \]
\end{Statement}

In order to prove the existence of sections of the Green-Griffiths logarithmic jet bundle, it is sufficient to prove the existence of sections of \(\mathcal{L}(a_{1},\dotsc,a_{\kappa})\) for a certain suitable choice of the parameters \(a_{1},\dotsc,a_{\kappa}\). 
More generally, for any ample line bundle \(\mathcal{A}\):
\[
  \dim H^0\bigl(
  \widebar{X}_{0},
  \mathcal{E}_{\kappa,\mu(\underline{a})}V_{0}^{\star}(\log D_{0})
  \otimes
  \mathcal{A}^{\vee}\bigr)
  \geq
  \dim H^0\bigl(
  \widebar{X}_{\kappa},
  \mathcal{L}(\underline{a})
  \otimes
  \pi_{\kappa,0}^{\star}\mathcal{A}^{\vee}
  \bigr).
\]

Starting with an hypersurface \(H\subset\P^{n}\) of degree \(d\), fix the log pair:
\[
  (\widebar{X}_{0},D_{0})
  \bydef
  (\P^{n},H).
\]
and also fix the distribution \(V_{0}\) by taking the whole logarithmic tangent space:
\[
  V_{0}
  \bydef
  T_{\P^{n}}(-\log H).
\]
Now, let as above denote the Demailly tower constructed on the log directed manifold \(\bigl(\P^{n},H,T_{\P^{n}}(-\log H)\bigr)\) by:
\[
  \bigl(\widebar{X}_{\kappa},D_{\kappa}\bigr)
  \to 
  \bigl(\widebar{X}_{\kappa-1},D_{\kappa-1}\bigr)
  \to \dotso\to
  \bigl(\widebar{X}_{1},D_{1}\bigr)
  \to
  \bigl(\widebar{X}_{0},D_{0}\bigr).
\] 
This context will be called the \textsl{logarithmic absolute case}.

Then using the positivity properties of \(\O_{X_{\kappa}(\underline{a}})\) and \(\O_{\P^{n}}(1)\), one can establish the following statement (\cite{MR2441250,MR2495771,MR1492539,MR1824906}):
\begin{Statement} 
  \label{prop:nef_cone}
  In the logarithmic absolute case, if \(a_{1},\dotsc,a_{\kappa}\) are \(\kappa\) positive integers
  having weighted sum
  \(
  \mu(\underline{a})
  =
  1\,a_{1}+\dotsb+\kappa\,a_{\kappa},
  \)
  and satisfying the inequalities:
  \[
    a_{i}
    \geq
    3\,a_{i+1}
    >0
    \pour{i=1,\dotsc,\kappa-1},
  \]
  then the line bundle:
  \[
    \mathcal{F}
    \bydef
    (\pi_{\kappa,0})^{\star}
    \O_{\P^{n}}\bigl(2\mu(\underline{a})\bigr)
    \otimes
    \mathcal{L}\bigl(a_{1},\dotsc,a_{\kappa}\bigr)
    \to
    \widebar{X}_{\kappa}
  \]
  is nef.
\end{Statement}
Recall~\cite[1.4.1]{MR2095471} that a line bundle \(\mathcal{L}\to X\) over \(X\) is numerically effective (\textsl{nef}) if for every irreducible curve \(\mathscr{C}\subset X\): 
\[
  \int_{\mathscr{C}}
  c_{1}(\mathcal{L})
  \geq
  0.
\]
The second factor \(\mathcal{L}(a_{1},\dotsc,a_{\kappa})\) of \(\mathcal{F}\) is positive on curves lying in the fibers, and the first factor \(\O_{\P^{n}}(2\mu(\underline{a})\) compensate the negativity than could occur in the horizontal directions.

\subsection{Morse inequalities}
Now, the line bundle at stake:
\[
  \mathcal{L}(a_{1},\dotsc,a_{\kappa})
  \otimes 
  (\pi_{\kappa,0})^{\star}
  \mathcal{K}_{\P^{n}}(H)^{-2\mu\delta}
\]
is the difference of two nef line bundles \(\mathcal{F},\mathcal{G}\) because one has the trivial identity (where the pullbacks are omitted):
\begin{multline*}
  \mathcal{L}(a_{1},\dotsc,a_{\kappa})
  \otimes 
  %(\pi_{\kappa,0})^{\star}
  \mathcal{K}_{\P^{n}}(H)^{-2\mu\delta}
  =\\
  \bigl(
  \underbrace{
    \mathcal{L}(a_{1},\dotsc,a_{\kappa})
    \otimes
    %(\pi_{\kappa,0})^{\star}
    \O_{\P^{n}}(2\mu)
  }_{\defby\mathcal{F}}
  \bigr)
  \otimes
  \bigl(
  \underbrace{
    %(\pi_{\kappa,0})^{\star}
    \O_{\P^{n}}\bigl(
    2\mu+2\mu\delta(d-n-1))
    \bigr)
  }_{\defby\mathcal{G}}
  \bigr)^{\vee}.
\end{multline*}
Indeed, it has already be said just above that \(\mathcal{F}\) is nef, and \(\mathcal{G}\) is also nef because it is the pullback of a nef line bundle.

Together with the following Demailly-Trapani \textsl{Morse inequalities}(\cite{MR1339712,MR1492539}) this fact provides with a (rough) control the cohomology of the line bundle \(\mathcal{L}(a_{1},\dotsc,a_{\kappa})\to\widebar{X}_{\kappa}\). 

\begin{Statement}
  \label{thm_morse}
  For any holomorphic line bundle \(\mathcal{L}\) on a \(N\)-dimensional compact manifold \(X\), that  can be written as the difference \(\mathcal{L}=\mathcal{F}\otimes\mathcal{G}^\vee\) of two nef line bundles \(\mathcal{F}\) and \(\mathcal{G}\), one has:
  \[
    \dim H^0\bigl(X,\mathcal{L}^{\otimes k}\bigr)
    \geq
    k^{N}\;
    \frac{
      \bigl(\mathcal{F}^{N}\bigr)-k\,\bigl(\mathcal{F}^{N-1}\cdot\mathcal{G}\bigr) 
    }
    {N!}
    - 
    o(k^{N}).
  \]
\end{Statement}
Recall that by definition~\cite{MR2095471}, the intersection number \(\bigl(L_{1}\dotsm L_k\bigr)\) of \(k\) line bundles on a \(k\) dimensional variety \(X\) denotes:
\[
  \bigl(L_{1}\dotsm L_k\bigr)
  \bydef
  \int_X
  c_{1}(L_{1})\dotsm c_{1}(L_k).
\]

These algebraic Morse inequalities \eqref{thm_morse} can now be applied to the constructed line bundle \(\mathcal{L}(\underline{a})\to\widebar{X}_{\kappa}\). By induction the dimension \(n_{\kappa}\bydef \dim(\widebar{X}_{\kappa})\) is simply computed as:
\[
  n_{\kappa}
  =
  \dim(\widebar{X}_{0})
  +
  \kappa\;\bigl(\mathrm{rk}\,P(V_{0})\bigr)
  =
  n+\kappa\,(n-1).
\]
\begin{Statement}
  If the integers \(a_{1},\dotsc,a_{\kappa}\) satisfy the inequalities:
  \[
    a_{i}
    \geq
    3\,a_{i+1}
    >
    0
    \pour{i=1,\dotsc,\kappa-1}
  \]
  and \(\delta\in\Q\) is a fixed positive rational number,
  It suffices to establish the positivity of the intersection number: 
  \[
    I
    \bydef
    \int_{\widebar{X}_{\kappa}}
    c_{1}(\mathcal{F})^{n_{\kappa}}
    -
    n_{\kappa}\,
    c_{1}(\mathcal{F})^{n_{\kappa}-1}c_{1}(\mathcal{G}),
  \]
  where:
  \[
    \left\{%}
    \begin{aligned}
      c_{1}(\mathcal{F})
      &=
      2\mu(\underline{a})\,h+a_{1}\,c_{1}(\mathcal{L}_{1}^{\vee})+\dotsb+a_{\kappa}\,c_{1}(\mathcal{L}_{\kappa}^{\vee})
      \\
      c_{1}(\mathcal{G})
      &=
      \bigl(1+\delta(d-n-1)\bigr)2\mu(\underline{a})\,h
    \end{aligned}
    \right.,
  \]
  in order to guarantee the existence of non zero global sections in the vector space:
  \[
    H^{0}\Bigl(
    \widebar{X}_{\kappa},
    \mathcal{E}_{\kappa,m}T_{\P^{n}}^{\star}\bigl(\log H\bigr)
    \otimes
    (\pi_{\kappa,0})^{\star}
    \bigl(\mathcal{K}_{\P^{n}}(H)^{\vee}\bigr)^{2\delta m}
    \Bigr),
  \] 
  for asymptotic weighted degree \(m\gg1\).
\end{Statement}
\begin{proof}
  Assume that one can find \(\underline{a}=(a_{1},\dotsc,a_{\kappa})\), such that \(I>0\) for \(d\geq\lambda(\underline{a})\).
  By the above algebraic Morse inequality, one has, for large integers \(k\):
  \[
    \dim H^0
    \Bigl(\widebar{X}_{\kappa},
    \bigl(
    \mathcal{L}(\underline{a})
    \otimes
    \mathcal{K}_{\P^{n}}(H)^{-2\delta\mu(\underline{a})}
    \bigr)^{\otimes k}
    \Bigr)
    \geq
    k^{n_{\kappa}}\;
    \frac{I}
    {n_{\kappa}!}
    - 
    o(k^{n_{\kappa}}).
  \]
  Now, both \(\mathcal{L}(\underline{a})\) and \(\mu(\underline{a})\) are linear in \(\underline{a}\), thus:
  \[
    \bigl(
    \mathcal{L}(\underline{a})
    \otimes
    \mathcal{K}_{\P^{n}}(H)^{-2\delta\mu(\underline{a})}
    \bigr)^{\otimes k}
    =
    \mathcal{L}(k\underline{a})
    \otimes
    \mathcal{K}_{\P^{n}}(H)^{-2\delta\mu(k\underline{a})},
  \]
  and \(k\underline{a}\) still satisfy the inequalities:
  \[
    k\,a_{i}
    \geq
    3\,k\,a_{i+1}
    >
    0
    \pour{i=1,\dotsc,\kappa-1}
  \]
  thus \(\mathcal{L}(k\underline{a})\) is a subsheaf of the Green-Griffiths logarithmic jet bundle:
  \[
    \mathcal{L}(k\underline{a})
    \subset
    \mathcal{E}_{\kappa,m}T_{\P^{n}}^{\star}(\log H),
  \]
  for the weighted degree:
  \[
    m
    \bydef
    \mu(k\underline{a})
    =
    k\mu(\underline{a}).
  \]

  By taking \(k\) large enough, one has:
  \[
    \dim H^0
    \Bigl(\widebar{X}_{\kappa},
    \mathcal{L}(k\underline{a})
    \otimes
    \mathcal{K}_{\P^{n}}(H)^{-2\delta\mu(k\underline{a})}
    \Bigr)
    \geq
    \frac{1}{2}\;
    k^{n_{\kappa}}\;
    \frac{I}
    {n_{\kappa}!}
    >0.
  \]
  Thus, for \(m\gg1\) large enough, \emph{not} effective:
  \[
    \dim H^0
    \Bigl(\widebar{X}_{\kappa},
    \mathcal{E}_{\kappa,m}T_{\P^{n}}^{\star}(\log H)
    \otimes
    \mathcal{K}_{\P^{n}}(H)^{-2\delta m}
    \Bigr)
    \geq
    1.
  \]
\end{proof}

Denote the first Chern classes of the vertical line bundles \(\mathcal{L}_{i}\) by:
\[
  v_{i}
  \bydef
  c_{1}\bigl(
  \mathcal{L}_{i}^{\vee}
  \bigr)
  \pour{i=1,\dotsc,\kappa}.
\]
then the integrand of the above intersection product:
\[
  f(v_{1},\dotsc,v_{\kappa})
  \bydef
  c_{1}(\mathcal{F})^{n_{\kappa}}
  -
  n_{\kappa}\,
  c_{1}(\mathcal{F})^{n_{\kappa}-1}c_{1}(\mathcal{G})
\]
can be written in this notation as:
\begin{multline*}
  f(v_{1},\dotsc,v_{\kappa})
  \bydef
  \Bigl(
  2\mu h+a_{1}v_{1}+\dotsb+a_{\kappa}v_{\kappa}
  \Bigr)^{n_{\kappa}}
  -
  \\
  n_{\kappa}\,
  \Bigl(
  \bigl(1+\delta(d-n-1)\bigr)
  2\mu h
  \Bigr)
  \Bigl(
  2\mu h+a_{1}v_{1}+\dotsb+a_{\kappa}v_{\kappa}
  \Bigr)^{n_{\kappa}-1}.
\end{multline*}
We will next give a formula for integrating \(f\) under this form.

\subsection{Formal computation on the Demailly tower}
It is convenient to bring down the computation to the basis. In \cite{D2013}, we provide a formula in that aim.
This formula expresses intersection products on \(\widebar{X}_{\kappa}\) as coefficients of an \textsl{iterated Laurent series}, in the very spirit of the predating residue formula of Berczi~(\cite{arxiv:1011.4710}). 

An iterated Laurent series is a multivariate formal series having well ordered support with respect to the lexicographic order on \(\Z^{n}\).
The main advantage of using the subspace of iterated Laurent series over using the whole space of multivariate formal series is that the Cauchy product is well defined.
There is an injection of the field of rational functions in the field of iterated Laurent series, called the iterated Laurent series expansion at the origin.
The reader is referred to \cite{D2013} for a precise description of this process.
By convention, in the product of an iterated Laurent series and a rational function, the rational function will always be replaced by its expansion. Accordingly, such a product is in fact the Cauchy product of two iterated Laurent series.
Lastly, the expression \([M]\varPhi\) will denote the coefficient of a monomial \(M\) in the multivariate formal series \(\varPhi\).

Applying the formula of \cite{D2013}, one can eliminate simultaneously all the vertical first Chern classes \(v_{i}\). 
\begin{Statement}
  \label{prop:residueFormula}
  In the ``logarithmic absolute case'', for any homogeneous polynomial: 
  \[
    f
    \in
    \C[h][t_{1},\dotsc,t_{\kappa}]
    =
    \C[h,t_{1},\dotsc,t_{\kappa}]
  \]
  having total degree \(n_{\kappa}=\dim\,\widebar{X}_{\kappa}\) (with respect to the variables \(h,t_{1},\dotsc,t_{\kappa}\)),
  the corresponding cohomology class
  \(
  f(v_{1},\dotsc,v_{\kappa})
  \in
  H^{\smbullet}(\widebar{X}_{\kappa})
  \)
  can be integrated on \(\widebar{X}_{\kappa}\) using the formula:
  \[
    \int_{\widebar{X}_{\kappa}}
    f(v_{1},\dotsc,v_{\kappa})
    =
    \bigl[
      h^{n}t_{1}^{n}\dotsm t_{\kappa}^{n}
    \bigr]
    \Bigl(
    \mathsf{A}(t_{1},\dotsc,t_{\kappa})\,
    \mathsf{B}(t_{1},\dotsc,t_{\kappa})\,
    \mathsf{C}(t_{1},\dotsc,t_{\kappa})
    \Bigr),
  \]
  where \(\mathsf{A}\) is the polynomial with coefficients in \(\C[h]\) obtained from \(f\) by:
  \[
    \mathsf{A}(t_{1},\dotsc,t_{\kappa})\,
    \bydef
    (dh+t_{1})
    \dotsm
    (dh+t_{\kappa})\;
    f(t_{1},\dotsc,t_{\kappa}),
  \]
  where \(\mathsf{B}\) is the Laurent polynomial with coefficients in \(\C[h]\):
  \[
    \mathsf{B}(t_{1},\dotsc,t_{\kappa})
    \bydef
    \sum_{j_{1},\dotsc,j_{\kappa}\geq0}
    \binom{n+j_{1}}{n}
    \dotsm
    \binom{n+j_{\kappa}}{n}\;
    \frac{(-h)^{j_{1}+\dotsb+j_{\kappa}}}{t_{1}^{j_{1}}\dotsm t_{\kappa}^{j_{\kappa}}},
  \]
  and where \(\mathsf{C}\) is the (unequivocal) iterated Laurent series expansion of the universal rational function:
  \[
    \mathsf{C}(t_{1},\dotsc,t_{\kappa})
    \bydef
    \prod_{1\leq i<j\leq\kappa}
    \left(
    \frac{t_j-t_i}{t_j-2t_i}
    \right)\;
    \prod_{2\leq i<j\leq\kappa}
    \left(
    \frac{t_j-2\,t_i}{t_j-2\,t_i+t_{i-1}}
    \right).
  \]
  with respect to the order: 
  \[
    t_{1}\ll\dotsb\ll t_{\kappa}\ll1.
  \] 
\end{Statement}
\begin{proof}
  This result is a plain implementation of the method of integration developed in~\cite{D2013}.
  In analogy with Chern polynomial, for a vector bundle \(E\to X\) over a \(N\) dimensional manifold \(X\), define \(s_u(E)\) to be the generating function of the Segre classes of \(E\), that is:
  \[
    s_u(E)
    \bydef
    s_{0}(E)+u\,s_{1}(E)+u^2\,s_{2}(E)+\dotsb+u^{n}\,s_N(E).
  \]
  In the precedent work \cite{D2013} we establish that:
  \begin{multline*}
    \int_{\widebar{X}_{\kappa}}
    f(v_{1},\dotsc,v_{\kappa})
    =
    \bigl[
      t_{1}^{n-1}\dotsm t_{\kappa}^{n-1}
    \bigr]
    \Biggl(
    \\
    \int_{\P^{n}}
    f(t_{1},\dotsc,t_{\kappa})\,
    \prod_{j=1}^{\kappa}
    s_{1/t_{j}}\bigl(T_{\P^{n}}^{\star}(\log D)\bigr)
    \\
    \prod_{1\leq i<j\leq\kappa}
    \left(
    \frac{t_j-t_i}{t_j-2t_i}
    \right)\;
    \prod_{2\leq i<j\leq\kappa}
    \left(
    \frac{t_j-2\,t_i}{t_j-2\,t_i+t_{i-1}}
    \right)
    \Biggr).
  \end{multline*}
  In the third line, we recognize the universal series \(\mathsf{C}\):
  \begin{multline*}
    \label{eq:fiberInt}
    \tag{\ensuremath{*}}
    \int_{\widebar{X}_{\kappa}}
    f(v_{1},\dotsc,v_{\kappa})
    =
    \\
    \bigl[
      t_{1}^{n-1}\dotsm t_{\kappa}^{n-1}
    \bigr]
    \Biggl(
    \int_{\P^{n}}
    f(t_{1},\dotsc,t_{\kappa})\,
    \prod_{j=1}^{\kappa}
    s_{1/t_{j}}\bigl(T_{\P^{n}}^{\star}(\log D)\bigr)\;
    \mathsf{C}(t_{1},\dotsc,t_{\kappa})
    \Biggr).
  \end{multline*}

  By a classic computation, the total Segre class of the basis is:
  \[
    s_{\smbullet}(V_{0})
    =
    s_{\smbullet}\bigl(T_{\P^{n}}(-\log H)\bigr)
    =
    \frac{(1+d\,h)}{(1+h)^{n+1}}.
  \]
  Thus the appearing Laurent polynomials \(s_{1/t_j}\bigl(T_{\P^{n}}^{\star}(\log D)\bigr)\) have the expression:
  \[
    s_{1/t_{i}}(V_{0})
    =
    \frac{t_{i}+dh}{t_{i}}\;
    \sum_{j_{i}=0}^{n}
    \binom{n+j_{i}}{j_{i}}
    \left(
    \frac{-h}{t_{i}}
    \right)^{\!\!j_{i}},
  \]
  and the product of these expressions for \(i=1,\dotsc,\kappa\) is:
  \[
    \prod_{i=1}^{\kappa}
    s_{1/t_{i}}\bigl(T_{\P^{n}}^{\star}(\log D)\bigr)\;
    =
    \frac{t_{1}+dh}{t_{1}}
    \dotsm
    \frac{t_{\kappa}+dh}{t_{\kappa}}\;
    \mathsf{B}(t_{1},\dotsc,t_{\kappa}).
  \]

  Then \eqref{eq:fiberInt} becomes:
  \begin{multline*}
    \int_{\widebar{X}_{\kappa}}
    f(v_{1},\dotsc,v_{\kappa})
    =
    \\
    \bigl[
      t_{1}^{n-1}\dotsm t_{\kappa}^{n-1}
    \bigr]
    \Biggl(
    \frac{1}{t_{1}\dotsm t_{\kappa}}
    \int_{\P^{n}}
    f(t_{1},\dotsc,t_{\kappa})\,
    \prod_{i=1}^{\kappa}
    (dh+t_{i})\;
    \mathsf{B}(t_{1},\dotsc,t_{\kappa})\;
    \mathsf{C}(t_{1},\dotsc,t_{\kappa})
    \Biggr).
  \end{multline*}
  We recognize the term \(\mathsf{A}\) and use that, obviously:
  \[
    \bigl[
      t_{1}^{n-1}\dotsm t_{\kappa}^{n-1}
    \bigr]
    \left(
    \frac{\varPhi}{t_{1}\dotsm t_{\kappa}}
    \right)
    =
    \bigl[
      t_{1}^{n}\dotsm t_{\kappa}^{n}
    \bigr]
    \Bigl(
    \varPhi
    \Bigr),
  \]
  in order to obtain:
  \[
    \int_{\widebar{X}_{\kappa}}
    f(v_{1},\dotsc,v_{\kappa})
    =
    \bigl[
      t_{1}^{n}\dotsm t_{\kappa}^{n}
    \bigr]
    \Biggl(
    \int_{\P^{n}}
    \mathsf{A}(t_{1},\dotsc,t_{\kappa})\;
    \mathsf{B}(t_{1},\dotsc,t_{\kappa})\;
    \mathsf{C}(t_{1},\dotsc,t_{\kappa})
    \Biggr).
  \]

  Now, it remains to integrate polynomials in the hyperplane class \(h\) on \(\P^{n}\). 
  For degree reason, the integrand is a multiple of \(h^{n}\). 
  The coefficient of this monomial is the sought intersection product, because it is known that:
  \[
    \int_{\P^{n}}h^{n}
    =
    1.
  \]

  Finally, our computation becomes the very concrete combinatorial problem:
  \[
    \int_{X_{\kappa}}
    f(v_{1},\dotsc,v_{\kappa})
    =
    \bigl[h^{n}t_{1}^{n}\dotsm t_{\kappa}^{n}\bigr]
    \Biggl(
    \mathsf{A}(t_{1},\dotsc,t_{\kappa})\;
    \mathsf{B}(t_{1},\dotsc,t_{\kappa})\;
    \mathsf{C}(t_{1},\dotsc,t_{\kappa})
    \Biggr).
  \]
\end{proof}

Moreover, this formula stay true without assumption on the degree of \(f\), because the only power of \(h\) that does not vanish by integration on the basis is the \(n\)-th power.

\section{Implementation of the computation}
\setcounter{equation}{0}
In this section, it is established that, if \(m\gg d\geq52n^{n}\) and if \(0<35n^{n}\delta\leq1\) then:
\[
  \dim H^0\Bigl(
  \mathcal{E}_{n,m}T_{\P^{n}}^{\star}\bigl(\log H\bigr)
  \otimes
  \mathcal{K}_{\P^{n}}(H)^{-2m\delta}
  \Bigr)
  \geq
  1,
\]
In that aim, it is sufficient to prove ---\,under the same hypotheses\,--- the positivity of the intersection product:
\[
  I
  =
  \bigl[h^{n}t_{1}^{n}\dotsm t_{\kappa}^{n}\bigr]
  \Bigl(
  \mathsf{A}(\underline{t})\,
  \mathsf{B}(\underline{t})\,
  \mathsf{C}(\underline{t})
  \Bigr),
\]
for a fixed choice of weights \(a_{1},a_{2},\dotsc,a_{\kappa}\), where \(\underline{t}\) stands for the \(\kappa\)-tuple \(t_{1},\dotsc,t_{\kappa}\).

In order to estimate the displayed Cauchy product coefficient \(I\), one needs to access to the coefficients of the terms \(\mathsf{A}\), \(\mathsf{B}\) and \(\mathsf{C}\).

\subsection{Preliminary expansion of \ensuremath{\mathsf{A}}}
The first appearing term is the polynomial:
\[
  \mathsf{A}(t_{1},\dotsc,t_{\kappa})
  =
  (d\,h+t_{1})
  \dotsm
  (d\,h+t_{\kappa})\
  f(t_{1},\dotsc,t_{\kappa}).
\]
It is the only term involving \(d\). Recall that in our situation, the polynomial \(f\) to be integrated on the Demailly tower is:
\begin{multline*}
  f\left(t_{1},\dotsc,t_{\kappa}\right)
  \bydef
  \Bigl(
  2\mu(\underline{a})h+a_{1}t_{1}+\dotsb+a_{\kappa}t_{\kappa}
  \Bigr)^{n_{\kappa}}
  -
  \\
  n_{\kappa}
  \Bigl(\delta\,d+(1-\delta\,n-\delta)\Bigr)
  \,
  2\mu(\underline{a})h
  \Bigl(
  2\mu(\underline{a})h+a_{1}t_{1}+\dotsb+a_{\kappa}t_{\kappa}
  \Bigr)^{n_{\kappa}-1},
\end{multline*}
where \(n_{\kappa}\) is the dimension of the ambient space at the \(\kappa\)-th level:
\[
  n_{\kappa}
  =
  \dim\widebar{X}_{\kappa}
  =
  n+\kappa(n-1),
\]
and where \(\mu\) is the weighted sum of the coefficients \(a_{i}\):
\[
  \mu(\underline{a})
  =
  1\,a_{1}+\dotsb+\kappa\,a_{\kappa}.
\]

Thus \(\mathsf{A}\) has degree at most \(n\) with respect to \(d\), because in its expression \(d\) appears only as a factor of the product \(dh\), and \(h^{n+1}=0\).

Let us introduce the short notation:
\[
  f_{i}(t_{1},\dotsc,t_{\kappa})
  \bydef
  \frac{n_{\kappa}!}{(n_{\kappa}-i)!}
  \bigl(2\mu(\underline{a})\bigr)^{i}
  \bigl(a_{1}t_{1}+\dotsc+a_{\kappa}t_{\kappa}\bigr)^{n_{\kappa}-i}
  \pour{i=0,1,\dotsc,n}
\]
for the terms (not depending on \(d\)) that appear in the expansion of \(f\) with respect to the hyperplane class \(h\):
\[
  f(t_{1},\dotsc,t_{\kappa})
  =
  \sum_{i=0}^{n}
  (\alpha_{i}-d\beta_{i})\,
  h^{i}\,
  f_{i}(t_{1},\dotsc,t_{\kappa}),
\]
where the rational coefficients \(\alpha_{i}\) and \(\beta_{i}\) have the respective expressions:
\[
  \alpha_{i}
  \bydef
  \frac{1-(1-\delta\,n-\delta)i}{i!}
  \quad\text{and}\quad
  \beta_{i}
  \bydef
  \frac{\delta i}{i!}
  \pour{i=0,1,\dotsc,n}.
\]
In particular, it will prove useful to compute that:
\[
  \alpha_{0} = 1,\;
  \alpha_{1} = \delta(n+1),\;
  \beta_{0} = 0,\;
  \beta_{1} = \beta_{2} = \delta.
\]

With this notation, one can write more shortly the expansion of the coefficients:
\[
  \mathsf{A}_{p}(t_{1},\dotsc,t_{\kappa})
  \bydef
  \bigl[d^{n-p}\bigr]
  \mathsf{A}(t_{1},\dotsc,t_{\kappa})
  \pour{p=0,1,\dotsc,n}.
\]
One has, by using the notation \(\mathfrak{s}_{j}\) for the \(j\)-th elementary symmetric function of \(\kappa\) variables: 
\[
  (d\,h+t_{1})
  \dotsm
  (d\,h+t_{\kappa})
  =
  \sum_{j=0}^{\kappa}
  (d\,h)^{\kappa-j}
  \mathfrak{s}_{j}(t_{1},\dotsc,t_{\kappa}).
\]
Thus, by definition, \(\mathsf{A}\) is the following product of polynomials in the variables \(d,h,t_{1},\dotsc,t_{\kappa}\):
\[
  \mathsf{A}(t_{1},\dotsc,t_{\kappa})
  =
  \sum_{i=0}^{n}
  \bigl(
  (\alpha_{i}-d\,\beta_{i})\,
  h^{i}\,
  f_{i}(\underline{t})
  \bigr)
  \sum_{j=\kappa-n}^{\kappa}
  \bigl(
  d^{\kappa-j}\,
  h^{\kappa-j}\,
  \mathfrak{s}_{j}(\underline{t})
  \bigr).
\]

Remind that \(\beta_{0}=0\), thus in this expression \(d\) appears \emph{always} as a factor of \(h\).
This remark yields at once:
\[
  \bigl[d^{n-p}\bigr]
  \mathsf{A}(t_{1},\dotsc,t_{\kappa})
  =
  O\bigl(h^{n-p}\bigr),
\]
where \(O\bigl(h^{n-p}\bigr)\) denotes a polynomial multiple of \(h^{n-p}\) in \(\C[h,t_{1},\dotsc,t_{\kappa}]\).

More precisely:
\begin{Statement}
  \label{eq:coeffA}
  For \(p=0,1,\dotsc,n\), the expansion of \(\mathsf{A}_{p}\) with respect to \(h\) is:
  \[
    \mathsf{A}_{p}
    =
    \sum_{q=0}^p
    h^{n-q}
    \Bigl(
    \alpha_{p-q}\,
    f_{p-q}(\underline{t})\;
    \mathfrak{s}_{\kappa-n+p}(\underline{t})
    -
    \beta_{p+1-q}\,
    f_{p+1-q}(\underline{t})\;
    \mathfrak{s}_{\kappa-n+p+1}(\underline{t})
    \Bigr),
  \]
  where \(\mathfrak{s}_{i}\) denotes the \(i\)-th elementary symmetric function of \(\kappa\) variables. 
\end{Statement}
\begin{proof}
  Recall:
  \[
    \mathsf{A}(t_{1},\dotsc,t_{\kappa})
    =
    \sum_{i=0}^{n}
    \bigl(
    (\alpha_{i}-d\,\beta_{i})\,
    h^{i}\,
    f_{i}(\underline{t})
    \bigr)
    \sum_{j=\kappa-n}^{\kappa}
    \bigl(
    d^{\kappa-j}\,
    h^{\kappa-j}\,
    \mathfrak{s}_{j}(\underline{t})
    \bigr).
  \]
  There is two ways to obtain \(d^{n-p}\); indeed the first polynomial is linear in \(d\). 
  Either take the constant coefficient (with respect to \(d\)) in the first polynomial and the coefficient of \(d^{n-p}\) in the second polynomial, that is:
  \[
    \tag{\(*\)}
    \sum_{i=0}^{p}
    \alpha_i\,
    h^i\,
    f_i(\underline{t})
    \cdot
    h^{n-p}\,
    \mathfrak{s}_{\kappa-n+p}(\underline{t}),
  \]
  or take the slope (with respect to \(d\)) in the first polynomial and the coefficient of \(d^{n-p-1}\) in the second polynomial, that is:
  \[
    \tag{\(**\)}
    \sum_{i=0}^{p+1}
    -\beta_{i}\,
    h^{i}\,
    f_{i}(\underline{t})
    \cdot
    h^{n-(p+1)}\,
    \mathfrak{s}_{\kappa-n+p+1}(\underline{t}).
  \]

  For any integer \(q=0,1,\dotsc,p\), 
  the coefficient of \(h^{n-q}\) in the first contribution (\(*\)) is:
  \[
    \alpha_{p-q}\;
    f_{p-q}(\underline{t})\;
    \mathfrak{s}_{\kappa-n+p}(\underline{t}),
  \]
  and the coefficient of \(h^{n-q}\) in the second contribution (\(**\)) is:
  \[
    -\beta_{p+1-q}\;
    f_{p+1-q}(\underline{t})\;
    \mathfrak{s}_{\kappa-n+p+1}(\underline{t}).
  \]
  One can easily deduce the expansion of \(\mathsf{A}_{p}\):
  \[
    \mathsf{A}_{p}
    =
    \sum_{q=0}^{p}
    h^{n-q}
    \Bigl(
    \alpha_{p-q}\,
    f_{p-q}(\underline{t})\,
    \mathfrak{s}_{\kappa-n+p}(\underline{t})
    -
    \beta_{p+1-q}\,
    f_{p+1-q}(\underline{t})\,
    \mathfrak{s}_{\kappa-n+p+1}(\underline{t})
    \Bigr).
  \]
  That is the announced result.
\end{proof}

\subsection{Examination of the other terms}
For multi-indices \(\mathbi{j}\in\Z^{\kappa}\), define the coefficients:
\[
  \mathsf{B}_{\mathbi{j}}
  \bydef
  (-1)^{j_{1}+\dotsb+j_{\kappa}}
  \binom{n+j_{1}}{n}
  \dotsm
  \binom{n+j_{\kappa}}{n}
\]
in such way that:
\[
  \mathsf{B}
  =
  \sum_{\mathbi{j}\in\N^{\kappa}}
  \mathsf{B}_{\mathbi{j}}\,
  \Biggl(\frac{h}{t_{1}}\Biggr)^{\!\!j_{1}}
  \dotsm
  \Biggl(\frac{h}{t_{\kappa}}\Biggr)^{\!\!j_{\kappa}}.
\]

The coefficients \(\mathsf{B}_{\mathbi{j}}\) have a simple shape.
In the contrary, it is highly challenging to understand the Laurent series expansion of the rational expression:
\[
  \mathsf{C}(t_{1},\dotsc,t_{\kappa})
  \bydef
  \prod_{1\leq i<j\leq\kappa}
  \left(
  \frac{t_{j}-t_{i}}{t_{j}-2t_{i}}
  \right)\;
  \prod_{2\leq i<j\leq\kappa}
  \left(
  \frac{t_{j}-2\,t_{i}}{t_{j}-2\,t_{i}+t_{i-1}}
  \right).
\]
In order to bypass this technical sticking point, in a first approximation, we will replace both terms \(\mathsf{B}\) and \(\mathsf{C}\) by their constant term \(1\).
Later, we will come back in more details to this point.

We use the notation:
\[
  \mathsf{C}_{\mathbi{k}}
  \bydef
  \bigl[
    t_{1}^{k_{1}}
    \dotsm
    t_{\kappa}^{k_{\kappa}}
  \bigr]
  \mathsf{C}(t_{1},\dotsc,t_{\kappa}).
\]
for the (complicated) coefficient of the monomial \(t_{1}^{k_{1}}\dotsm t_{\kappa}^{k_{\kappa}}\) in the iterated Laurent series expansion of \(\mathsf{C}\). By homogeneity, these coefficients are \(0\) unless \(k_{1}+\dotsb+k_{\kappa}=0\), in such way that:
\[
  \mathsf{C}(t_{1},\dotsc,t_{\kappa})
  =
  \sum_{k_{1}+\dotsb+k_{\kappa}=0}
  \mathsf{C}_{\mathbi{k}}\,
  t_{1}^{k_{1}}
  \dotsm
  t_{\kappa}^{k_{\kappa}}.
\]
Notice that one can use the rational expression of \(\mathsf{C}\) whenever \(\underline{t}\) is in the domain of convergence of \(\mathsf{C}\).
We will also use the transparent notation:
\[
  \abs{\mathsf{C}}
  =
  \sum_{k_{1}+\dotsb+k_{\kappa}=0}
  \abs{\mathsf{C}_{\mathbi{k}}}\,
  t_{1}^{k_{1}}
  \dotsm
  t_{\kappa}^{k_{\kappa}}.
\]

The shape of \(\mathsf{C}\) allows to extract some more information on the support of its iterated Laurent series expansion:
\begin{Statement}
  \label{lem:supportC}
  The coefficient \(\mathsf{C}_{\mathbi{k}}\) is zero unless for each \(i=1,\dotsc,\kappa\):
  \[
    k_{i}+\dotsb+k_{\kappa}
    \leq
    0.
  \]
\end{Statement}
\begin{proof}
  At each step of the iterated Laurent series expansion of \(\mathsf{C}\), one manipulates only products of the monomials: \(t_{i}/t_{j}\) for indices \(1\leq i<j\leq\kappa\).
\end{proof}

\subsection{Guide line of the computation}
For the moment, \emph{trust} that ---\,in a first approximation\,--- both terms \(\mathsf{B}\) and \(\mathsf{C}\) can be replaced by \(1\).
Take also \(\delta=0\). We introduce the auxiliary (positive) constants:
\begin{equation}
  \label{eq:auxiliaryIp}
  \widetilde{I}_{p}
  \bydef
  \bigl[t_{1}^{n}\dotsm t_{\kappa}^{n}\bigr]
  \Bigl(
  f_{p}(t_{1},\dotsc,t_{\kappa})\,
  \mathfrak{s}_{\kappa-n+p}(t_{1},\dotsc,t_{\kappa})
  \Bigr)
  \pour{p=0,1,\dotsc,n}.
\end{equation}
Compare with the actual coefficients of \(d^{n-p}\) in \(I\):
\begin{small}
  \[
    \widetilde{I}_{p}
    \bydef
    \bigl[h^{n}t_{1}^{n}\dotsm t_{\kappa}^{n}\bigr]
    \left(
    \sum_{q=0}^{p}
    h^{n-q}
    \Bigl(
    \alpha_{p-q}\,
    f_{p-q}(\underline{t})\,
    \mathfrak{s}_{\kappa-n+p}(\underline{t})
    -
    \beta_{p+1-q}\,
    f_{p+1-q}(\underline{t})\,
    \mathfrak{s}_{\kappa-n+p+1}(\underline{t})
    \Bigr)
    \right).
  \]
\end{small}

\medskip

Although this seems to be a radical simplification, It will be proven later that for all \(p=0,1,\dotsc,n\), the coefficient \(I_{p}\) has nearly the same size as the simplified coefficient \(\widetilde{I}_{p}\), provided the parameters \(a_{1},\dotsc,a_{\kappa}\) are \emph{suitably adjusted}~\eqref{hyp:1} and that \(\delta\ll1\) is taken \emph{small enough}~\eqref{hyp:2}.

\smallskip
Here is the guide line of our computation:
\par\noindent\(\square\) 
\emph{Firstly}, we will compute a sufficient lower bound \(\widetilde{\lambda}(\underline{a})\) on the degree \(d\), such that the approximated polynomial:
\[
  \widetilde{I}(d)
  \bydef
  \widetilde{I}_{0}\,d^{n}-\bigl(\widetilde{I}_{1}\,d^{n-1}+\dotsb+\widetilde{I}_{n}\bigr)
\]
has positive values for degrees \(d\geq\widetilde{\lambda}(\underline{a})\).

\begin{small}
  \noindent\sffamily
  The goal of all the technicalities that will follow is to bring the computation back to this \emph{much more} easy computation.
\end{small}
\smallskip

\noindent\(\square\)
\emph{Secondly}, we will study the leading coefficient, and prove that under \textit{ad hoc} hypotheses \eqref{hyp:1} and \eqref{hyp:2}, it has a positive value, that is at least two thirds of \(\widetilde{I}_{0}\):
\smallskip\newline\noindent\frame{
  \parbox \textwidth{
    \begin{multline*}
      I_{0}
      =
      \underbrace{
        \bigl[
          h^{n}t_{1}^{n}\dotsm t_{\kappa}^{n}
        \bigr]
        \Bigl(
        \mathsf{A}_{0}(t_{1},\dotsc,t_{\kappa})\vert_{\delta=0}
        \Bigr)
      }_{=\widetilde{I}_{0}}
      +\\
      \underbrace{
        \bigl[
          h^{n}t_{1}^{n}\dotsm t_{\kappa}^{n}
        \bigr]
        \Bigl(
        \mathsf{A}_{0}(t_{1},\dotsc,t_{\kappa})\vert_{\delta=0}\,
        \bigl(\mathsf{C}(t_{1},\dotsc,t_{\kappa})-1\bigr)
        \Bigr)
      }_{\text{negligible for \(a_{1}\gg\dotsb\gg a_{\kappa}\)}\;\eqref{hyp:1}}
      +
      \underbrace{
        \vphantom{\Big(}
        \quad
        \delta\,I_{0}'
        \quad
      }_{\text{negligible for \(\delta\ll1\)}\;\eqref{hyp:2}}.
    \end{multline*}
  }
}
\smallskip

\noindent\(\square\)
\emph{Thirdly}, we will study the remaining coefficients and derive a sufficient lower bound \(\lambda(\underline{a})\) on the degree \(d\) ---\,similar to \(\widetilde{\lambda}(\underline{a})\)\,--- from the relative sizes of the coefficients.
Again, we will see that:
\smallskip\newline\noindent\frame{
  \parbox \textwidth{
    \begin{multline*}
      I_{p}
      =
      \underbrace{
        \bigl[
          h^{n}t_{1}^{n}\dotsm t_{\kappa}^{n}
        \bigr]
        \Bigl(
        \mathsf{A}_{p}(\underline{t})\vert_{\delta=0}\,
        \mathsf{B}(\underline{t})
        \Bigr)
      }_{\sim\text{multiple of \(\widetilde{I}_{p}\)}}
      +\\
      \underbrace{
        \bigl[
          h^{n}t_{1}^{n}\dotsm t_{\kappa}^{n}
        \bigr]
        \Bigl(
        \mathsf{A}_{p}(\underline{t})\vert_{\delta=0}\,
        \mathsf{B}(\underline{t})\,
        \bigl(\mathsf{C}(\underline{t})-1\bigr)
        \Bigr)
      }_{\text{negligible for \(a_{1}\gg\dotsb\gg a_{\kappa}\)}\;\eqref{hyp:1}}
      +
      \underbrace{
        \vphantom{\Big(}
        \quad
        \delta\,I_{p}'
        \quad
      }_{\text{negligible for \(\delta\ll1\)}\;\eqref{hyp:2}}.
    \end{multline*}
  }
}

\subsection{Quick justification of the approximations}
Recall that \(I_{p}\) is the Cauchy product coefficient:
\[
  I_{p}
  =
  \bigl[
    h^{n}t_{1}^{n}\dotsm t_{\kappa}^{n}
  \bigr]
  \Bigl(
  \mathsf{A}_{p}\left(\underline{t}\right)\,
  \mathsf{B}\left(\underline{t}\right)\,
  \mathsf{C}\left(\underline{t}\right)
  \Bigr)
\]
where \(\mathsf{A}_{p}\) is the polynomial:
\[
  \mathsf{A}_{p}
  =
  \sum_{q=0}^{p}
  h^{n-q}
  \Bigl(
  \alpha_{p-q}\,
  f_{p-q}(\underline{t})\,
  \mathfrak{s}_{\kappa-n+p}(\underline{t})
  -
  \beta_{p+1-q}\,
  f_{p+1-q}(\underline{t})\,
  \mathfrak{s}_{\kappa-n+p+1}(\underline{t})
  \Bigr),
\]
and \(\mathsf{B}\) and \(\mathsf{C}\) are the iterated Laurent series:
\[
  \mathsf{B}
  =
  \sum_{\mathbi{j}}
  \mathsf{B}_{\mathbi{j}}\,
  \Biggl(\frac{h}{t_{1}}\biggr)^{\!\!j_{1}}
  \dotsm
  \Biggl(\frac{h}{t_{\kappa}}\Biggr)^{\!\!j_{\kappa}}
  \qquad\text{and:}\qquad
  \mathsf{C}
  =
  \sum_{\mathbi{k}}
  \mathsf{C}_{\mathbi{k}}\,
  t_{1}^{k_{1}}\dotsm t_{\kappa}^{k_{\kappa}}.
\]
As a consequence, \(I_{p}\) can be written as:
\begin{multline*}
  I_{p}
  =
  \sum_{\mathbi{i}=\mathbi{j}-\mathbi{k}}\;
  \sum_{0\leq q\leq p}
  \alpha_{p-q}\,
  \bigl[
    t_{1}^{n+i_{1}}\dotsm t_{\kappa}^{n+i_{\kappa}}
  \bigr]
  \Bigl(
  f_{p-q}(\underline{t})\,
  \mathfrak{s}_{\kappa-n+p}(\underline{t})
  \Bigr)\,
  \mathsf{B}_{\mathbi{j}}\,
  \mathsf{C}_{\mathbi{k}}
  -\\
  \sum_{\mathbi{i}=\mathbi{j}-\mathbi{k}}\;
  \sum_{0\leq q\leq p}
  \beta_{p+1-q}\,
  \bigl[
    t_{1}^{n+i_{1}}\dotsm t_{\kappa}^{n+i_{\kappa}}
  \bigr]
  \Bigl(
  f_{p+1-q}(\underline{t})\,
  \mathfrak{s}_{\kappa-n+p+1}(\underline{t})
  \Bigr)\,
  \mathsf{B}_{\mathbi{j}}\,
  \mathsf{C}_{\mathbi{k}}.
\end{multline*}

In order to make the estimation of this term easier, the only freedom that we have is the choice of the parameters 
\(a_{1}>\dotsb>a_{\kappa}>0\) 
that appear in:
\[
  f_{i}(t_{1},\dotsc,t_{\kappa})
  =
  \bigl(2\mu(\underline{a})\bigr)^{i}
  \frac{n_{\kappa}!}{(n_{\kappa}-i)!}
  \bigl(a_{1}t_{1}+\dotsb+a_{\kappa}t_{\kappa}\bigr)^{n_{\kappa}-i}.
\]
This latitude will become more clear after we explain how to simplify the appearing coefficients:
\[
  \bigl[
    t_{1}^{n+i_{1}}\dotsm t_{\kappa}^{n+i_{\kappa}}
  \bigr]
  \Bigl(
  f_{p-q}(t_{1},\dotsc,t_{\kappa})\;
  \mathfrak{s}_{\kappa-n+p}(t_{1},\dotsc,t_{\kappa})
  \Bigr),
\]
in order to deal with more tractable plain \(\underline{a}\)-monomials. 
%In particular, the term \(\mathfrak{s}_{\kappa-n}(\underline{t})\) prove to be a slight technical complication.
Notice that for degree reasons, such coefficients are zero unless:
\[
  i_{1}+\dotsb+i_{\kappa}
  =
  q
\]

Remind the auxiliary constant introduced in \eqref{eq:auxiliaryIp}:
\[
  \widetilde{I}_{p}
  =
  \bigl[
    t_{1}^{n}\dotsm t_{\kappa}^{n}
  \bigr]
  \Bigl(
  f_{p}(t_{1},\dotsc,t_{\kappa})\,
  \mathfrak{s}_{\kappa-n+p}(t_{1},\dotsc,t_{\kappa})\;
  \Bigr).
\]
It stands for the absolute value of the coefficients of \(d^{n-p}\) in the simplified polynomial \(\widetilde{I}\), appearing in the preceding guide line. Each coefficient appearing in the computation of the coefficient of \(d^{n-p}\) in \(I\) can be compared with the unique coefficient \(\widetilde{I}_{p}\) appearing in the computation of the coefficient of \(d^{n-p}\) in \(\widetilde{I}\):
\begin{Statement}
  \label{lem:fp}
  For any integers \(p,q\in\N\) and any multi-index 
  \(\bigl(i_{1},\dotsc,i_{\kappa}\bigr)\in\Z^{\kappa}\), 
  such that 
  \[
    i_{1}+\dotsb+i_{\kappa}
    =
    q
    \leq
    p,
  \]
  one has:
  \[
    \bigl[
      t_{1}^{n+i_{1}}\dotsm t_{\kappa}^{n+i_{\kappa}}
    \bigr]
    \Bigl(
    f_{p-q}(\underline{t})\;
    \mathfrak{s}_{\kappa-n+p}(\underline{t})
    \Bigr)
    \leq
    \left(
    \frac{a_{1}}{2n\mu(\underline{a})}
    \right)^{i_{1}}
    \dotsm
    \left(
    \frac{a_{\kappa}}{2n\mu(\underline{a})}
    \right)^{i_{\kappa}}\;
    \widetilde{I}_{p}.
  \]
\end{Statement}
\begin{proof}
  By definition, the \((\kappa-n+p)\)-th elementary symmetric function is the following homogeneous sum of monomials of total degree \(\kappa-n+p\):
  \[
    \tag{\(*\)}
    \mathfrak{s}_{\kappa-n+p}(t_{1},\dotsc,t_{\kappa})
    =
    \sum_{\substack{(\varepsilon_{1},\dotsc,\varepsilon_{\kappa})\in\{0,1\}^\kappa\\\varepsilon_{1}+\dotsb+\varepsilon_{\kappa}=\kappa-n+p}}
    t_{1}^{\varepsilon_{1}}\dotsm t_{\kappa}^{\varepsilon_{\kappa}}.
  \]
  Remind that \(q=(i_{1}+\dotsb+i_{\kappa})\) and \(n_{\kappa}=\kappa\,n-(\kappa-n)\). Accordingly, the degree of the monomial that we consider:
  \[
    (n+i_{1})+\dotsb+(n+i_{\kappa})
    =
    \kappa\,n+q
  \] 
  agrees with he homogeneous degree of the expression at stake:
  \[
    \bigl(n_{\kappa}-(p-q)\bigr)+(\kappa-n+p)
    =
    \kappa\,n-(\kappa-n+p)+q+(\kappa-n+p).
  \]
  By definition:
  \[
    f_{j}(\underline{t})
    =
    \frac{n_{\kappa}!\,(2\mu)^{j}}{(n_{\kappa}-j)!}
    (a_{1}t_{1}+\dotsb+a_{\kappa}t_{\kappa})^{n_{\kappa}-j}.
  \]
  This expression can be expanded thanks to the multinomial formula, and by replacing \(f_{p-q}\) by its expansion, one gets the formula:
  \[
    &\bigl[
      t_{1}^{n+i_{1}}\dotsm t_{\kappa}^{n+i_{\kappa}}
    \bigr]
    \Bigl(
    f_{p-q}(\underline{t})
    \mathfrak{s}_{\kappa-n+p}(\underline{t})
    \Bigr)
    \\
    &\quad=
    \sum_{\substack{
      (\varepsilon_{1},\dotsc,\varepsilon_{\kappa})\in\{0,1\}^{\kappa}\\
      \varepsilon_{1}+\dotsb+\varepsilon_{\kappa}=\kappa-n+p
    }}
    \frac{n_{\kappa}!\;(2\mu)^{p-q}}{(n_{\kappa}-p+q)!}
    \frac{(n_{\kappa}-p+q)!}{(n+i_{1}-\varepsilon_{1})!\dotsm(n+i_{\kappa}-\varepsilon_{\kappa})!}
    a_{1}^{n+i_{1}-\varepsilon_{1}}\dotsm a_{\kappa}^{n+i_{\kappa}-\varepsilon_{\kappa}},
    \\
    &\quad=
    \sum_{\substack{(\varepsilon_{1},\dotsc,\varepsilon_{\kappa})\in\{0,1\}^\kappa\\\varepsilon_{1}+\dotsb+\varepsilon_{\kappa}=\kappa-n+p}}
    \frac{n_{\kappa}!\;(2\mu)^{p-q}}{(n+i_{1}-\varepsilon_{1})!\dotsm(n+i_{\kappa}-\varepsilon_{\kappa})!}
    a_{1}^{n+i_{1}-\varepsilon_{1}}\dotsm a_{\kappa}^{n+i_{\kappa}-\varepsilon_{\kappa}},
  \]
  where it should be understood that we consider only the \(\kappa\)-tuples \((\varepsilon_{1},\dotsc,\varepsilon_{\kappa})\) such that for \(k=1,\dotsc,\kappa\) the integer \(n+i_{k}-\varepsilon_{k}\) is not negative. 

  By distinguishing the case \((i<0)\) and the case \((i>0)\) it is easy to show that in both case:
  \[
    \frac{n!}{(n+i)!}
    \leq
    n^{-i}.
  \]

  Thus, by summing the powers of \(n\):
  \[
    \frac{1}{(n+i_{1}-\varepsilon_{1})!\dotsm(n+i_{\kappa}-\varepsilon_{\kappa})!}
    \leq
    \frac{n^{\kappa-n+p+q}}{(n!)^{\kappa}}.
  \]
  In the right hand side, the dependence on \(\varepsilon_{1},\dotsc,\varepsilon_{\kappa}\) disappears. The monomials can thus be summed by using \((*)\) in the other direction, and one gets:
  \begin{multline*}
    \bigl[
      t_{1}^{n+i_{1}}\dotsm t_{\kappa}^{n+i_{\kappa}}
    \bigr]
    \Bigl(
    f_{p-q}(t_{1},\dotsc,t_{\kappa})
    \mathfrak{s}_{\kappa-n+p}(t_{1},\dotsc,t_{\kappa})
    \Bigr)
    \\\leq
    \frac{n^{\kappa-n+p-q}\;n_{\kappa}!}{(n!)^{\kappa}}
    (2\mu)^{p-q}a_{1}^{n+i_{1}}\dotsm a_{\kappa}^{n+i_{\kappa}}
    \mathfrak{s}_{\kappa-n+p}(1/a_{1},\dotsc,1/a_{\kappa}).
  \end{multline*}
  It remains only to compare with the expression of the auxiliary constant \(\widetilde{I}_{p}\) in order to conclude.
  In intermediary computation, its explicit value was computed to be:
  \[
    \widetilde{I}_{p}
    =
    \frac{n^{\kappa-n+p}\,n_{\kappa}!}{(n!)^{\kappa}}\;
    (2\mu)^p\;
    (a_{1}\dotsm a_{\kappa})^{n}
    \mathfrak{s}_{\kappa-n+p}(1/a_{1},\dotsc,1/a_{\kappa}).
  \]
  Thus:
  \[
    \bigl[
      t_{1}^{n-i_{1}}\dotsm t_{\kappa}^{n-i_{\kappa}}
    \bigr]
    \Bigl(
    f_{p-q}(\underline{t})
    \mathfrak{s}_{\kappa-n+p}(\underline{t})
    \Bigr)
    \leq
    \widetilde{I}_{p}\;
    \frac{a_{1}^{i_{1}}\dotsm a_{\kappa}^{i_{\kappa}}}
    {(2n\mu)^{q}}.
  \]
\end{proof}

We will use this latitude and set \(a_{1}\gg\dotsb\gg a_{\kappa}\) in order to make the annoying contributions
negligible with respect to the ``central'' term \(\widetilde{I}_{p}\), obtained for \(\mathbi{i}=(0,\dotsc,0)\in\Z^{\kappa}\).

Consider the coefficients:
\[
  \left(
  \frac{a_{1}}{2n\mu(\underline{a})}
  \right)^{i_{1}}
  \dotsm
  \left(
  \frac{a_{\kappa}}{2n\mu(\underline{a})}
  \right)^{i_{\kappa}},
\]
and use the telescoping products: 
\[
  \left(\frac{a_{j}}{2n\mu(\underline{a})}\right)
  =
  \left(\frac{a_{j}}{a_{j-1}}\right)
  \dotsm
  \left(\frac{a_{2}}{a_{1}}\right)
  \left(\frac{a_{1}}{2n\mu(\underline{a})}\right)
\]
in order to write it the coefficients under the form:
\[
  \left(
  \frac{a_{1}}{2n\mu(\underline{a})}
  \right)^{i_{1}+\dotsb+i_{\kappa}}
  \left(
  \frac{a_{2}}{a_{1}}
  \right)^{i_{2}+\dotsb+i_{\kappa}}
  \dotsm
  \left(
  \frac{a_{\kappa}}{a_{\kappa-1}}
  \right)^{i_{\kappa}}.
\]
The coefficient of \(\mathbi{i}=(0,\dotsc,0)\) will always be \(1\), but for other indices, for parameters \(a_{i}\gg a_{i+1}>0\) ---\,that clearly fulfil the condition \eqref{prop:nef_cone}\,---, these coefficients can be made as small as wanted.
Indeed:
\begin{Statement}
  all multi-indices \(\mathbi{i}\) involved in this computation have entries such that:
  \[
    i_{j}+\dotsb+i_{\kappa}
    \geq
    0
    \pour{j=1,\dotsc,\kappa}.
  \]
\end{Statement}
\begin{proof}
  One has \(\mathbi{i}=\mathbi{j}-\mathbi{k}\). The entries of \(\mathbi{j}\) have to be non negative and we have seen that \(\mathsf{C}_{\mathbi{k}}\) is zero unless:
  \[
    k_{i}+\dotsb+k_{\kappa}
    \leq
    0
    \pour{i=1,\dotsc,\kappa}.
  \]
\end{proof}

Our understanding does not allow to give the geometric meaning of this hypothesis.
The careful study of preceding works \cite{MR2593279,arxiv:1011.4710} let appear that this kind of hypothesis is also done in order to get a result in any dimension.
We do not know if there is a similar hypothesis in \cite{arxiv:1011.3636}.
For all that, this technical simplification yields the desired statement.

\emph{However}, if one is interested in effective computation, the parameters \(a_{1},\dotsc,a_{\kappa}\) have to be chosen with care. Indeed, a very quickly decreasing sequence \(\underline{a}\) will guarantee a very large leading coefficient, but will also increase the size of the largest positive root of the polynomial \(I\).

\subsection{Simplified computation}
Now, we achieve our program. As a first step, let us simplify the problem and let us consider the polynomial inequality:
\[
  \widetilde{I}_{0}\,d^{n}
  \geq
  \widetilde{I}_{1}\,d^{n-1}+\dotsb+\widetilde{I}_{n-1}\,d+\widetilde{I}_{n}.
\]
Because \(\widetilde{I}_{0}>0\), this inequality is checked as soon as \(d\) is larger than the (largest) positive root of the polynomial:
\[
  \widetilde{I}(d)
  =
  \widetilde{I}_{0}\,d^{n}-\widetilde{I}_{1}\,d^{n-1}-\dotsb-\widetilde{I}_{n-1}\,d-\widetilde{I}_{n}.
\]
Dividing the polynomial \(\widetilde{I}\) by the appearing (huge) positive constant \(\widetilde{I}_{n}\), later computed to be:
\[
  \widetilde{I}_{n}
  =
  (2\mu)^{n}\,
  \frac{n_{\kappa}!}{(n-1)!^\kappa}\,
  (a_{1}\dotsm a_{\kappa})^{n-1},
\]
does not change the values for which \(\widetilde{I}\) is positive.  
In order to bound from above the absolute values of the roots of a polynomial,
it suffices to estimate the relative size of its coefficients, according to the following lemma:
\begin{Statement}
  All complex roots of a non constant complex algebraic equation of degree \(n\):
  \[
    c_{0}x^{n}=c_{1}x^{n-1}+\dotsb+c_{n-1}x+c_n
  \]
  have an absolute value least than or equal to the \textsl{Fujiwara's bound}:
  \[
    \lambda
    =
    2\,\max_{p=1,\dotsc,n}
    \bigg\lvert
    \frac{c_p}{c_{0}}
    \bigg\rvert^{1/p}.
  \]
\end{Statement}

Following this strategy of bounding \(\widetilde{I}_{p}/\widetilde{I}_{0}\), we will now compute the relative sizes of the coefficients of \(\widetilde{I}\) and after that compare each non-leading coefficient with the leading coefficient:
\begin{Statement}
  \label{lem:cp}
  For any integer \(p=0,1,\dotsc,n\), one has:
  \[
    \widetilde{I}_{p}
    =
    \mathfrak{s}_{n-p}
    \left(
    \frac{a_{1}}{2n\mu},
    \dotsc,
    \frac{a_{\kappa}}{2n\mu}
    \right)
    \widetilde{I}_{n}.
  \]
\end{Statement}
\begin{proof}
  By definition:
  \[
    f_{p}(t_{1},\dotsc,t_{\kappa})
    =
    (2\mu)^{p}\,
    \frac{n_{\kappa}!}{(n_{\kappa}-p)!}\,
    \bigl(a_{1}\,t_{1}+\dotsb+a_{\kappa}\,t_{\kappa}\bigr)^{n_{\kappa}-p}.
  \]
  We will expand this \((n_{\kappa}-p)\)-power. In this view, recall that for any two integers \(k,l\in\N\) the standard binomial formula yield by induction on \(k\) the following \textsl{multinomial formula}:
  \[
    \bigl(x_{1}+\dotsb+x_{k}\bigr)^l
    =
    \sum_{i_{1}+\dotsb+i_{k}=l}
    \frac{l!}{i_{1}!\dotsm i_{k}!}\,
    x_{1}^{i_{1}}\dotsm x_{k}^{i_{k}}.
  \]
  Thus \(f_{p}\) becomes:
  \[
    f_{p}(t_{1},\dotsc,t_{\kappa})
    =
    (2\mu)^{p}\,
    \sum_{i_{1}+\dotsb+i_{\kappa}=n_{\kappa}-p}
    \frac{n_{\kappa}!}{(n_{\kappa}-p)!}
    \frac{(n_{\kappa}-p)!}{i_{1}!\dotsm i_{\kappa}!}\;
    a_{1}^{i_{1}}\dotsm a_{\kappa}^{i_{\kappa}}\;
    t_{1}^{i_{1}}\dotsm t_{\kappa}^{i_{\kappa}}.
  \]

  On the other hand, by definition, the \((\kappa-n+p)\)-th elementary symmetric function of \(\kappa\) variables \(x_{1},\dotsc,x_{`}\) is the following homogeneous sum of monomials: 
  \[
    \tag{\(*\)}
    \mathfrak{s}_{\kappa-n+p}(x_{1},\dotsc,x_{\kappa})
    =
    \sum_{\substack{
      (\varepsilon_{1},\dotsc,\varepsilon_{\kappa})\in\{0,1\}^{\kappa}\\
      \varepsilon_{1}+\dotsb+\varepsilon_{\kappa}=\kappa-n+p
    }}
    x_{1}^{\varepsilon_{1}}
    \dotsm 
    x_{\kappa}^{\varepsilon_{\kappa}},
  \]
  the total degree of which is \(\kappa-n=p\) and
  in which the exponent of each variable \(x_{1},\dotsc,x_{\kappa}\) is at most \(1\).

  The product of the two polynomials \(f_{p}(\underline{t})\) and \(\mathfrak{s}_{\kappa-n+p}(\underline{t})\) is thus:
  \[
    f_{p}(\underline{t})\,
    \mathfrak{s}_{\kappa-n+p}(\underline{t})
    =
    (2\mu)^{p}
    \sum_{\substack{
      i_{1}+\dotsb+i_{\kappa}=n_{\kappa}-p\\
      \varepsilon_{1}+\dotsb+\varepsilon_{\kappa}=\kappa-n+p
    }}
    \frac{n_{\kappa}!}{i_{1}!\dotsm i_{\kappa}!}\;
    a_{1}^{i_{1}}\dotsm a_{\kappa}^{i_{\kappa}}\;
    t_{1}^{i_{1}+\varepsilon_{1}}
    \dotsm 
    t_{\kappa}^{i_{\kappa}+\varepsilon_{\kappa}},
  \]
  where as above:
  \[
    (i_{1},\dotsc,i_{\kappa})
    \in
    \N^{\kappa},\;
    (\varepsilon_{1},\dotsc,\varepsilon_{\kappa})
    \in
    \{0,1\}^{\kappa}.
  \]
  Note that the homogeneous degree of the polynomial at stake agrees with the total degree of the monomial \(t_{1}^{n}\dotsm t_{\kappa}^{n}\) because it is:
  \[
    n_{\kappa}-p+(\kappa-n+p)
    =
    n+\kappa\,(n-1)-p+\kappa-n+p
    =
    \kappa\,n.
  \]

  Now, for \(i_{1}+\varepsilon_{1}=n,\dotsc,i_{\kappa}+\varepsilon_{\kappa}=n\), all \(i_{j}\) are non negative and the appearing multinomial coefficient has always the same value:
  \[
    \frac{n_{\kappa}!}{i_{1}!\dotsm i_{\kappa}!}
    =
    \frac{n_{\kappa}!}{(n-\varepsilon_{1})!\dotsm(n-\varepsilon_{\kappa})!}
    =
    \frac{n_{\kappa}!}{n!^{n-p}\;(n-1)!^{\kappa-n+p}}
    =
    \frac{n_{\kappa}!}{n^{n-p}\;(n-1)!^{\kappa}},
  \]
  because \(\kappa-n+p\) of the \(\varepsilon_{i}\)'s are \(1\) and the \(n-p\) remaining \(\varepsilon_{i}\)'s are \(0\).

  As a consequence:
  \[
    \bigl[t_{1}^{n}\dotsm t_{\kappa}^{n}\bigr]
    \Bigl(
    f_{p}(\underline{t})\,
    \mathfrak{s}_{\kappa-n+p}(\underline{t})
    \Bigr)
    =
    (2\mu)^{p}\!\!\!\!
    \sum_{\varepsilon_{1}+\dotsb+\varepsilon_{\kappa}=\kappa-n+p}
    \frac{n_{\kappa}!}{n^{n-p}\,(n-1)!^{\kappa}}
    a_{1}^{n-\varepsilon_{1}}\dotsm a_{\kappa}^{n-\varepsilon_{\kappa}}.
  \]

  Next, one can factorize the multinomial coefficient and \(a_{1}^{n}\dotsm a_{\kappa}^{n}\) and use the above definition (\(*\)) of the elementary symmetric function from the right to the left, in order to obtain:
  \[
    &
    \bigl[t_{1}^{n}\dotsm t_{\kappa}^{n}\bigr]
    \Bigl(
    f_{p}(\underline{t})\,
    \mathfrak{s}_{\kappa-n+p}(\underline{t})
    \Bigr)
    \\&\qquad
    =
    \frac{(2\mu)^{p}\;n_{\kappa}!}{n^{n-p}\,(n-1)!^{\kappa}}\;
    a_{1}^{n}\dotsm a_{\kappa}^{n}
    \sum_{\varepsilon_{1}+\dotsb+\varepsilon_{\kappa}=\kappa-n+p}
    \biggl(\frac{1}{a_{1}}\biggr)^{\varepsilon_{1}}
    \dotsm
    \biggl(\frac{1}{a_{\kappa}}\biggr)^{\varepsilon_{\kappa}}
    \\&\qquad
    =
    \frac{(2\mu)^{p}}{n^{n-p}}
    \frac{n_{\kappa}!}{(n-1)!^{\kappa}}\;
    a_{1}^{n}\dotsm a_{\kappa}^{n}\;\;
    \mathfrak{s}_{\kappa-n+p}\biggl(\frac{1}{a_{1}},\dotsc,\frac{1}{a_{\kappa}}\biggr).
  \]
  Finally we use that, by putting all fractions to the same denominator:
  \[
    \mathfrak{s}_{\kappa-n+p}
    \biggl(\frac{1}{a_{1}},\dotsc,\frac{1}{a_{\kappa}}\biggr)
    =
    \frac{1}{a_{1}\dotsm a_{\kappa}}\;
    \mathfrak{s}_{n-p}(a_{1},\dotsc,a_{\kappa}),
  \]
  and by homogeneity:
  \[
    \frac{(2\mu)^{p}}{n^{n-p}}\;
    \mathfrak{s}_{n-p}(a_{1},\dotsc,a_{\kappa})
    =
    (2\mu)^{n}\;
    \mathfrak{s}_{n-p}\biggl(\frac{a_{1}}{2n\mu},\dotsc,\frac{a_{\kappa}}{2n\mu}\biggr),
  \]
  hence we get:
  \[
    \widetilde{I}_{p}
    =
    (2\mu)^{n}\,
    \frac{n_{\kappa}!}{(n-1)!^{\kappa}}\,
    (a_{1}\dotsm a_{\kappa})^{n-1}\,
    \mathfrak{s}_{n-p}\biggl(\frac{a_{1}}{2n\mu},\dotsc,\frac{a_{\kappa}}{2n\mu}\biggr),
  \]
  that is the announced expression, because now:
  \[
    \widetilde{I}_{n}
    =
    (2\mu)^{n}\,
    \frac{n_{\kappa}!}{(n-1)!^{\kappa}}\,
    (a_{1}\dotsm a_{\kappa})^{n-1}.
  \]
  Thus:
  \[
    \widetilde{I}_{p}
    =
    \widetilde{I}_{n}\;
    \mathfrak{s}_{n-p}\biggl(\frac{a_{1}}{2n\mu},\dotsc,\frac{a_{\kappa}}{2n\mu}\biggr).
  \]
\end{proof}

As a consequence:
\begin{Statement}
  \label{cor:approximatedFujiwaraBound}
  The Fujiwara's bound for the polynomial \(\widetilde{I}\):
  \[
    \widetilde{\lambda}(a_{1},\dotsc,a_{\kappa})
    \bydef
    2\,\max_{1\leq p\leq n}
    \bigg\lvert
    \frac{\widetilde{I}_{p}}{\widetilde{I}_{0}}
    \bigg\rvert^{1/p}
  \]
  has the value:
  \[
    \widetilde{\lambda}(\underline{a})
    =
    4n\,
    \mu(\underline{a})\,
    \frac{\mathfrak{s}_{n-1}(\underline{a})}
    {\mathfrak{s}_{n}(\underline{a})}.
  \]
\end{Statement}
\begin{proof}
  By the lemma just above:
  \[
    \bigg\lvert
    \frac{\widetilde{I}_{p}}{\widetilde{I}_{0}}
    \bigg\rvert^{1/p}
    =
    \biggl(
    \frac{\mathfrak{s}_{n-p}(a_i/(2n\mu))}{\mathfrak{s}_n(a_i/(2n\mu))}
    \biggr)^{1/p}.
  \]

  \textsf{Claim:}
  For \(\kappa\in\N\) and \(p=1,\dotsc,\kappa-1\), and a set of positive reals \(x_{1},\dotsc,x_{\kappa}\):
  \[
    \mathfrak{s}_p(x_j)^2
    \geq
    \mathfrak{s}_{p-1}(x_j)\,\mathfrak{s}_{p+1}(x_j).
  \]
  In order to prove this claim, we will construct an injection from the set parametrizing \(\mathfrak{s}_{p+1}\,\mathfrak{s}_{p-1}\) to the set parametrizing \(\mathfrak{s}_p^2\); Let \(i_{1}<\dotsb<i_{p+1}\vert j_{1}<\dotsb<j_{p-1}\) be an element of the first set.
  \begin{itemize}
    \item 
      Either for any \(k=1,\dotsc,p-1\) one has \(i_k\geq j_k\). Then \(i_{p+1}>i_{p-1}\geq j_{p-1}\) and:
      \[
        i_{1}<\dotsb<i_p\vert j_{1}<\dotsb<j_{p-1}<i_{p+1}
      \]
      is in the second set.
    \item 
      Or there is a minimum \(k\) such that \(i_k<j_k\). Then if \(k>1\), \(j_{k-1}\leq i_{k-1}<i_k\) and in all cases \(i_k<j_k\), thus:
      \[
        i_{1}<\dotsb<i_{k-1}<i_{k+1}<\dotsb<i_p\vert j_{1}<\dotsb<j_{k-1}<i_k<j_k<\dotsb<j_{p-1}
      \]
      is in the second set.
  \end{itemize}
  The claim is now proved.
  It yields:
  \[
    \frac{\mathfrak{s}_{n-k}(x_j)}{\mathfrak{s}_{n-k+1}(x_j)}
    \leq
    \frac{\mathfrak{s}_{n-k+1}(x_j)}{\mathfrak{s}_{n-k+2}(x_j)}
    \leq
    \dotsb
    \leq
    \frac{\mathfrak{s}_{n-1}(x_j)}{\mathfrak{s}_{n}(x_j)},
  \]
  and thus:
  \[
    \biggl(
    \frac{\mathfrak{s}_{n-p}(x_j)}{\mathfrak{s}_n(x_j)}
    \biggr)^{\!1/p}
    =
    \biggl(
    \prod_{k=1}^p
    \frac{\mathfrak{s}_{n-k}(x_j)}{\mathfrak{s}_{n-k+1}(x_j)}
    \biggr)^{\!1/p}
    \leq
    \frac{\mathfrak{s}_{n-1}(x_j)}{\mathfrak{s}_{n}(x_j)}.
  \]
  One can directly deduce from this that:
  \[
    \lambda(a_{1},\dotsc,a_{\kappa})
    =
    2\,\max_{1\leq p\leq n}
    \biggl(
    \frac{\mathfrak{s}_{n-p}(a_i/(2n\mu))}
    {\mathfrak{s}_{n}(a_i/(2n\mu))}
    \biggr)^{\!1/p}
    =
    2\,
    \frac{\mathfrak{s}_{n-1}(a_{i}/(2n\mu))}
    {\mathfrak{s}_{n}(a_{i}/(2n\mu))}.
  \]
  Finally, by homogeneity:
  \[
    \widetilde{\lambda}(\underline{a})
    =
    4\,n\mu\,\frac{\mathfrak{s}_{n-1}(\underline{a})}{\mathfrak{s}_{n}(a_{j})}.
  \]
  This ends the computation.
\end{proof}

Using the above Fujiwara's bound, and taking account that the leading coefficient is positive, we obtain that the approximated  polynomial \(\widetilde{I}\) is positive for degrees:
\begin{equation}
  d
  \geq
  4\,n\;
  \mu(\underline{a})\;
  \frac{\mathfrak{s}_{n-1}(\underline{a})}
  {\mathfrak{s}_{n}(\underline{a})}.
\end{equation}

\bigskip
Next we will explain how to follow the same strategy with the (more complicated) actual computation.

\subsection{Estimation of the leading coefficient}
Actually, the leading coefficient \(I_{0}\) of the polynomial:
\[
  I(d)
  =
  I_{0}\,d^{n}-I_{1}\,d^{n-1}
  -\dotsb-
  I_{n-1}\,d-I_{n}
\]
is the Cauchy product coefficient:
\[
  I_{0}
  =
  \bigl[h^{n}t_{1}^{n}\dotsm t_{\kappa}^{n}\bigr]
  \Bigl(
  \mathsf{A}_{0}(t_{1},\dotsm,t_{\kappa})\,
  \mathsf{B}(t_{1},\dotsm,t_{\kappa})\,
  \mathsf{C}(t_{1},\dotsm,t_{\kappa})
  \Bigr),
\]
and we want to prove that it is \emph{positive} and to justify that it is comparable to the simpler coefficient:
\[
  \widetilde{I}_{0}
  =
  \bigl[t_{1}^{n}\dotsm t_{\kappa}^{n}\bigr]
  \Bigl(
  f_{0}(t_{1},\dotsc,t_{\kappa})\,
  \mathfrak{s}_{\kappa-n}(t_{1},\dotsc,t_{\kappa})
  \Bigr).
\]

There are indeed some similarities between the computations of these two coefficients (\textit{cf. supra} for the earlier computation of \(\widetilde{I}_{0}\)). 

Recall that \(\mathsf{A}_{0}\) is the coefficient of \(d^{n}\) in \(\mathsf{A}\) and it has the expression:
\[
  \mathsf{A}_{0}(t_{1},\dotsc,t_{\kappa})
  =
  h^{n}
  \Bigl(
  f_{0}\bigl(\underline{t}\bigr)\mathfrak{s}_{\kappa-n}\bigl(\underline{t}\bigr)
  -
  \delta\,f_{1}\bigl(\underline{t}\bigr)\mathfrak{s}_{\kappa-n+1}\bigl(\underline{t}\bigr)
  \Bigr),
\]
Thus \([h^{n}]\mathsf{A}_{0}(t_{1},\dotsc,t_{\kappa})\) and \(f_{0}(t_{1},\dotsc,t_{\kappa})\,\mathfrak{s}_{\kappa-n}(t_{1},\dotsc,t_{\kappa})\) coincide for \(\delta=0\).

Moreover, in this multivariate polynomial \(\mathsf{A}\), the exponent of \(h\) is already the same as the exponent of \(h\) in the sought monomial \(h^{n}t_{1}^{n}\dotsm t_{\kappa}^{n}\). Thus, the other term depending on \(h\):
\[
  \mathsf{B}(t_{1},\dotsc,t_{\kappa})
  =
  \sum_{i_{1},\dotsc,i_{\kappa}\geq0}
  (-1)^{i_{1}+\dotsb+i_{\kappa}}
  \binom{n+i_{1}}{n}
  \dotsm
  \binom{n+i_{\kappa}}{n}\,
  \frac{h^{i_{1}+\dotsb+i_{\kappa}}}
  {t_{1}^{i_{1}}\dotsm t_{\kappa}^{i_{\kappa}}},
\]
can be replaced by its truncation \(\mathsf{B}(\underline{t})=1+O(h)\) ---\,where \(O(h)\) denotes a polynomial multiple of \(h\) in \(\C[h,t_{1},\dotsc,t_{\kappa}]\)\,---, and the resulting simplified computation is:
\[
  I_{0}
  =
  \bigl[h^{n}t_{1}^{n}\dotsm t_{\kappa}^{n}\bigr]
  \Bigl(
  \mathsf{A}_{0}(h,t_{1},\dotsm,t_{\kappa})\,
  \mathsf{C}(t_{1},\dotsm,t_{\kappa})
  \Bigr),
\]
that is the sum:
\[
  I_{0}
  &=
  \sum_{k_{1}+\dotsb+k_{\kappa}=0}
  \bigl[h^{n}t_{1}^{n-k_{1}}\dotsm t_{\kappa}^{n-k_{\kappa}}\bigr]
  \Bigl(
  \mathsf{A}_{0}(h,t_{1},\dotsm,t_{\kappa})\,
  \Bigr)
  \bigl[t_{1}^{k_{1}}\dotsm t_{\kappa}^{k_{\kappa}}\bigr]
  \Bigl(
  \mathsf{C}(t_{1},\dotsm,t_{\kappa})
  \Bigr).
\]
That can be written as the difference of two sums:
\begin{multline*}
  I_{0}=
  \sum_{k_{1}+\dotsb+k_{\kappa}=0}
  \bigl[t_{1}^{n-k_{1}}\dotsm t_{\kappa}^{n-k_{\kappa}}\bigr]
  \Bigl(
  f_{0}\bigl(\underline{t}\bigr)\mathfrak{s}_{\kappa-n}\bigl(\underline{t}\bigr)
  \Bigr)
  \bigl[t_{1}^{k_{1}}\dotsm t_{\kappa}^{k_{\kappa}}\bigr]
  \Bigl(
  \mathsf{C}\bigl(\underline{t}\bigr)
  \Bigr)
  \\-\delta\,
  \sum_{k_{1}+\dotsb+k_{\kappa}=0}
  \bigl[t_{1}^{n-k_{1}}\dotsm t_{\kappa}^{n-k_{\kappa}}\bigr]
  \Bigl(
  f_{1}\bigl(\underline{t}\bigr)\mathfrak{s}_{\kappa-n+1}\bigl(\underline{t}\bigr)
  \Bigr)
  \bigl[t_{1}^{k_{1}}\dotsm t_{\kappa}^{k_{\kappa}}\bigr]
  \Bigl(
  \mathsf{C}\bigl(\underline{t}\bigr)
  \Bigr).
\end{multline*}
Of course, by taking \(\delta\ll1\) small enough, the second of these sums become negligible. 

It remains to study the first displayed sum. Here, there is a delicate pairing between the numerous appearing coefficients of the polynomial \(f_{0}\,\mathfrak{s}_{\kappa-n}\) and the corresponding coefficients of \(\mathsf{C}\). The combinatorics of the first family of coefficients is well understood, whereas the combinatorics of the second family of coefficients is \emph{very intricate}.

Notice that the coefficient \(\widetilde{I}_{0}\) corresponds to the single term of this sum indexed by \(\mathbi{i}=(0,\dotsc,0)\).
Thus, comparing \(I_{0}\) and \(\widetilde{I}_{0}\) amounts to replace \(\mathsf{C}\) by its constant term \(1\) (and also take \(\delta=0\)).

The only freedom that we have, in order to facilitate the estimation of that sum, is the choice of the parameters \(a_{1},\dotsc,a_{\kappa}\) that are involved in the terms \(f_{0}\) and \(f_{1}\) appearing in the expression:
\[
  \mathsf{A}_{0}(h,t_{1},\dotsc,t_{\kappa})
  =
  h^{n}
  \Bigl(
  f_{0}\bigl(\underline{t}\bigr)\mathfrak{s}_{\kappa-n}\bigl(\underline{t}\bigr)
  -
  \delta\,f_{1}\bigl(\underline{t}\bigr)\mathfrak{s}_{\kappa-n+1}\bigl(\underline{t}\bigr)
  \Bigr).
\]
We will use this latitude in the following proposition, in order to make all terms appearing in the sum:
\[
  \sum_{\substack{
    k_{1}+\dotsb+k_{\kappa}=0\\
    \mathbi{i}\neq(0,\dotsc,0)
  }}
  \bigl[t_{1}^{n-k_{1}}\dotsm t_{\kappa}^{n-k_{\kappa}}\bigr]
  \Bigl(
  f_{0}\bigl(\underline{t}\bigr)\mathfrak{s}_{\kappa-n}\bigl(\underline{t}\bigr)
  \Bigr)
  \bigl[t_{1}^{k_{1}}\dotsm t_{\kappa}^{k_{\kappa}}\bigr]
  \Bigl(
  \mathsf{C}\bigl(\underline{t}\bigr)
  \Bigr)
\]
negligible with respect to:
\[
  \widetilde{I}_{0}
  =
  \bigl[t_{1}^{n-0}\dotsm t_{\kappa}^{n-0}\bigr]
  \Bigl(
  f_{0}\bigl(\underline{t}\bigr)\mathfrak{s}_{\kappa-n}\bigl(\underline{t}\bigr)
  \Bigr)
  \bigl[t_{1}^{0}\dotsm t_{\kappa}^{i_{0}}\bigr]
  \Bigl(
  \mathsf{C}\bigl(\underline{t}\bigr)
  \Bigr).
\]

Below in appendix \ref{apx:leadingCoeff} we show the following statement:
\begin{Statement}
  \label{prop:leadingCoefficient}
  In the logarithmic absolute case, where \(V=T_{\P^{n}}(-\log H)\), fix the parameter \(\underline{a}\) by taking:

  \parbox{\textwidth}{
    \[
      \label{hyp:1}
      \tag{\ensuremath{\mathcal{H}_{1}}}
      \kappa=n\geq6
      \text{ and }
      \underline{a}=\bigl(n^{n},n^{n-1},\dotsc,n,1\bigr)
    \]
  }
  and \textit{a posteriori} take \(\delta=\delta(\underline{a})\) small enough:

  \parbox{\textwidth}{
    \[
      \label{hyp:2}
      \tag{\ensuremath{\mathcal{H}_{2}}}
      5\,\delta(\underline{a})\,\widetilde{\lambda}(\underline{a})
      \leq
      1,
    \]
  }
  then the leading coefficient of the polynomial \(I(d)\) is positive and it is at least two thirds of the coefficient \(\widetilde{I}_{0}\):
  \[
    I_{0}
    \geq
    \frac{2}{3}\,\widetilde{I}_{0}
    >0.
  \]
\end{Statement}

\subsection{Computation of the Fujiwara's bound}
Next, in order to compute the Fujiwara's bound of the polynomial \(I(d)\), we estimate the other coefficients. Considering the exponent \(1/p\) in the definition of the Fujiwara's bound, one is easily convinced that it is more important to have a good estimate on the coefficient \(I_{1}\) of the monomial \(d^{n-1}\). Thus we treat it separately.
Below in appendix \ref{apx:dn-1} we show the following statement:
\begin{Statement}
  \label{lem:boundI1}
  For any \(\kappa\geq n\) and any \(a_{1},\dotsc,a_{\kappa}\),
  the absolute value of the coefficient of \(d^{n-1}\) in the polynomial \(I\) is bounded from above by:
  \[
    \abs{I_{1}}
    \leq
    \widetilde{I}_{1}\;
    \left(
    \frac{1+3\,\delta\widetilde{\lambda}}{2}
    +
    \frac{1+\delta\widetilde{\lambda}}{2n}
    +
    \delta(n+1)
    \right)\;
    \abs{\mathsf{C}}
    \left(
    \frac{1}{a_{1}},
    \dotsc,
    \frac{1}{a_{\kappa}}
    \right).
  \]
  In particular, under hypotheses \eqref{hyp:1} and \eqref{hyp:2}:
  \[
    \abs{I_{1}}
    \leq
    5\,\widetilde{I}_{1}.
  \]
\end{Statement}

For the remaining coefficients, we do not expand the term \(\mathsf{B}\), thus we get a slightly larger upper bound.
Below in appendix \ref{apx:Ip} we show the following statement:
\begin{Statement}
  \label{lem:boundIp}
  For any \(\kappa\geq n\), and any \(a_{1},\dotsc,a_{\kappa}\), for \(p=0,1,\dotsc,n\), the absolute value of the coefficient of \(d^{n-p}\) in the polynomial \(I\) is bounded from above by:
  \[
    \abs{I_{p}}
    \leq
    \widetilde{I}_{p}\;
    \left(
    \frac{2+\delta\widetilde{\lambda}}{2}\,
    \left(
    \frac{2\,n}{2\,n-1}
    \right)^{{\!\!}n+1}
    \right)\;
    \abs{\mathsf{C}}
    \left(
    \frac{1}{a_{1}}
    ,\dotsc,
    \frac{1}{a_{\kappa}}
    \right).
  \]
  In particular, under hypotheses \eqref{hyp:1} and \eqref{hyp:2}:
  \[
    \abs{I_{p}}
    \leq
    12\,\widetilde{I}_{p}.
  \]
\end{Statement}

One can now compare the actual Fujiwara's bound \(\lambda\) with the approximated Fujiwara's bound \(\widetilde{\lambda}\). One has:
\begin{Statement}
  Under hypotheses \eqref{hyp:1} and \eqref{hyp:2}, the Fujiwara's bound:
  \[
    \lambda(\underline{a})
    \bydef
    2\,\max_{p=1,\dotsc,n}\left\lvert\frac{I_p}{I_{0}}\right\rvert^{1/p}
  \]
  of the polynomial \(I(d)\) satisfy:
  \[
    \lambda(\underline{a})
    \leq
    \frac{15}{2}\,\widetilde{\lambda}(\underline{a}).
  \]
\end{Statement}
\begin{proof}
  Recall that, by definition, \(\widetilde{\lambda}\) is the Fujiwara's bound of the simplified polynomial \(\widetilde{I}\):
  \[
    \widetilde{\lambda}(\underline{a})
    =
    2\,\max_{1\leq p\leq n}\left(\frac{\widetilde{I}_{p}}{\widetilde{I}_{0}}\right)^{\!\!1/p}.
  \]
  In \eqref{prop:leadingCoefficient} above, we have proven that, under hypotheses \eqref{hyp:1} and \eqref{hyp:2}, one has:
  \[
    I_{0}
    \geq
    \frac{2}{3}\,\widetilde{I}_{0}
  \]
  and in \eqref{lem:boundI1} just above that, under the same hypotheses:
  \[
    \abs{I_{1}}
    \leq
    5\,\widetilde{I}_{1}.
  \]
  Thus, by combining these two inequalities, we obtain:
  \[
    \frac{\abs{I_{1}}}{\abs{I_{0}}}
    \leq
    \frac{15}{2}
    \frac{\widetilde{I}_{1}}{\widetilde{I}_{0}}
  \]

  Similarly, in \eqref{lem:boundIp} we show that for \(p=2,\dotsc,n\):
  \[
    \abs{I_{p}}
    \leq
    12\,\widetilde{I}_{p}
  \]
  thus one has the corresponding inequality:
  \[
    \left(\frac{\abs{I_{p}}}{\abs{I_{0}}}\right)^{\!\!1/p}
    &\leq
    18^{1/p}\;
    \left(\frac{\widetilde{I}_{p}}{\widetilde{I}_{0}}\right)^{1/p}
    \leq
    3\sqrt{2}\;
    \left(\frac{\widetilde{I}_{p}}{\widetilde{I}_{0}}\right)^{1/p}.
  \]
\end{proof}

In order to receive a sufficient lower bound \(d_{n}\) on the degree \(d\), such that \(I(d)\) is positive for \(d\geq d_{n}\), it remains to estimate the value of \(\widetilde{\lambda}\) under our two hypotheses.
\begin{Statement}
  Under hypotheses \eqref{hyp:1} and \eqref{hyp:2}, the approximated Fujiwara's bound \(\widetilde{\lambda}(a_{1},a_{2},\dotsc,a_{n})\) satisfy:
  \[
    4
    \leq
    \frac{\widetilde{\lambda}(a_{1},\dotsc,a_{\kappa})}{n^{n}}
    \leq
    4\,\biggl(\frac{n}{n-1}\biggr)^{\!\!3}.
  \]
\end{Statement}
\begin{proof}
  Recall that:
  \[
    \widetilde{\lambda}(\underline{a})
    &=
    4\,n\mu(\underline{a})\,
    \frac{\mathfrak{s}_{n-1}(a_{i})}{\mathfrak{s}_{n}(a_{i})}.
  \]
  Now, \(n\mu(\underline{a})/n^{n}\) is the partial sum of the convergent power series with positive coefficients:
  \[
    \sum_{i\geq1}i\,(1/n)^{i-1},
  \]

  thus one has:
  \[
    n^{n}\cdot1
    \leq
    \left(
    n\mu(\underline{a})
    =
    n^{n}
    \cdot
    \sum_{i=1}^{n}
    i\left(\frac{1}{n}\right)^{i-1}
    \right)
    \leq
    n^{n}\cdot
    \left(\frac{n}{n-1}\right)^{\!\!2},
  \]
  and similarly:
  \[
    1
    \leq
    \left(
    \frac{\mathfrak{s}_{n-1}(a_{i})}{\mathfrak{s}_{n}(a_{i})}
    =
    \sum_{i=0}^{n-1}\biggl(\frac{1}{n}\biggr)^{i}
    \right)
    \leq
    \frac{n}{n-1}.
  \]
  These inequalities yield the stated bounds on \(\widetilde{\lambda}\).
\end{proof}

Then, we simplify the bounds on \(d\) and \(\delta\)  by using that for \(n\geq6\):
\[
  \begin{cases}
    \lambda
    \leq
    15/2*4*\left(\frac{n}{n-1}\right)^{\!\!3}
    n^{n}
    \leq
    52\,n^{n}
    \\
    1/\delta
    \geq
    35\,n^{n}
    \geq
    5*4*
    \left(\frac{n}{n-1}\right)^{\!\!3}
    n^{n}
    \Rightarrow
    5\delta\widetilde{\lambda}
    \leq
    1
  \end{cases}.
\]

\newpage
\appendix
\small
\section{Positivity of the leading coefficient}
\label{apx:leadingCoeff}
\begin{itshape}
  In this appendix, we establish how the technical hypotheses \eqref{hyp:1} and \eqref{hyp:2} yield the positivity of the leading coefficient \(I_{0}\) of \(I(d)\). We fix \(\kappa=n\).
\end{itshape}
\medskip
\subsection{Proof of proposition \eqref{prop:leadingCoefficient}}

The leading coefficient of \(I(d)\) has the expression:
\begin{multline*}
  \tag{\(*\)}
  I_{0}
  =
  \sum_{k_{1}+\dotsb+k_{n}=0}
  \bigl[t_{1}^{n-k_{1}}\dotsm t_{n}^{n-k_{n}}\bigr]
  \Bigl(
  f_{0}\bigl(\underline{t}\bigr)
  \Bigr)
  \bigl[t_{1}^{k_{1}}\dotsm t_{n}^{k_{n}}\bigr]
  \Bigl(
  \mathsf{C}\bigl(\underline{t}\bigr)
  \Bigr)
  \\-\delta\,
  \sum_{k_{1}+\dotsb+k_{n}=0}
  \bigl[t_{1}^{n-k_{1}}\dotsm t_{n}^{n-k_{n}}\bigr]
  \Bigl(
  f_{1}\bigl(\underline{t}\bigr)\mathfrak{s}_{1}\bigl(\underline{t}\bigr)
  \Bigr)
  \bigl[t_{1}^{k_{1}}\dotsm t_{n}^{k_{n}}\bigr]
  \Bigl(
  \mathsf{C}\bigl(\underline{t}\bigr)
  \Bigr).
\end{multline*}
We introduce some notation for the positive and negative contributions to the first sum, as follows:
\[
  I_{0}^{+}
  &\bydef
  \sum_{\substack{
    k_{1}+\dotsb+k_{n}=0\\
    [\underline{t}^{\mathbi{k}}\,]\mathsf{C}(\underline{t})>0
  }}
  \bigl[t_{1}^{n-k_{1}}\dotsm t_{n}^{n-k_{n}}\bigr]
  \Bigl(
  f_{0}\bigl(\underline{t}\bigr)
  \Bigr)
  \bigl[t_{1}^{k_{1}}\dotsm t_{n}^{k_{n}}\bigr]
  \Bigl(
  \mathsf{C}\bigl(\underline{t}\bigr)
  \Bigr),
  \\
  I_{0}^{-}
  &\bydef
  -
  \sum_{\substack{
    k_{1}+\dotsb+k_{n}=0\\
    [\underline{t}^{\mathbi{k}}\,]\mathsf{C}(\underline{t})<0
  }}
  \bigl[t_{1}^{n-k_{1}}\dotsm t_{n}^{n-k_{n}}\bigr]
  \Bigl(
  f_{0}\bigl(\underline{t}\bigr)
  \Bigr)
  \bigl[t_{1}^{k_{1}}\dotsm t_{n}^{k_{n}}\bigr]
  \Bigl(
  \mathsf{C}\bigl(\underline{t}\bigr)
  \Bigr),
\]
and also for the slope of \(I_{0}\) with respect to \(\delta\):
\[
  I_{0}'
  \bydef
  \sum_{k_{1}+\dotsb+k_{n}=0}
  \bigl[t_{1}^{n-k_{1}}\dotsm t_{n}^{n-k_{n}}\bigr]
  \Bigl(
  f_{1}(\underline{t})\,
  \mathfrak{s}_{1}(\underline{t})
  \Bigr)
  \bigl[t_{1}^{k_{1}}\dotsm t_{n}^{k_{n}}\bigr]
  \Bigl(
  \mathsf{C}\bigl(\underline{t}\bigr)
  \Bigr).
\]
Then the above equation \((*)\) becomes:
\[
  I_{0}
  =
  I_{0}^{+}-I_{0}^{-}-\delta\,I_{0}'.
\]
\smallskip

We prove below in sections \eqref{apx:C}, \eqref{apx:C_bound}, \eqref{apx:delta}, that under the hypothesis \eqref{hyp:1}, one has the following estimates.
\begin{itemize}
  \item 
    \emph{Firstly}, the sum of positive contributions is bounded from below by:
    \[
      I_{0}^{+}
      \geq
      2\,\widetilde{I}_{0}.
    \]
  \item
    \emph{Secondly}, the sum of negative contributions is bounded from above by:
    \[
      I_{0}^{-}
      \leq
      \bigl(5/6\bigr)\,\widetilde{I}_{0}.
    \]
  \item
    \emph{Thirdly}, the slope with respect to \(\delta\) is bounded from above by:
    \[
      \abs{I_{0}'}
      \leq
      \bigl(5\widetilde{\lambda}(\underline{a})/2\bigr)\;
      \widetilde{I}_{0}.
    \]
\end{itemize}
Putting these three estimates together, finishes the proof of \eqref{prop:leadingCoefficient}. Indeed, one obtains:
\[
  I_{0}
  \geq
  I_{0}^{+}
  -
  I_{0}^{-}
  -
  \delta\,
  \abs{I_{0}'}
  \geq
  \widetilde{I}_{0}
  \left(
  2
  -
  \frac{5}{6}
  -
  \frac{5\,\delta\widetilde{\lambda}}{2}
  \right)
  \geq
  \widetilde{I}_{0}
  \left(
  \frac{7-15\,\delta\widetilde{\lambda}}{6}
  \right).
\]
It suffices to take, as in \eqref{hyp:2}:
\[
  5\,
  \delta\bigl(\underline{a}\bigr)
  \widetilde{\lambda}\bigl(\underline{a}\bigr)
  \leq
  1,
\] 
in order to obtain as announced:
\[
  I_{0}
  \geq
  \frac{2}{3}\,\widetilde{I}_{0}
  >
  0.
\]

\subsection{Partial expansion of \ensuremath{\mathsf{C}} and positive contributions}
\label{apx:C}
In \eqref{lem:supportC}, we have seen that the coefficient \(\mathsf{C}_{\mathbi{k}}\) is zero unless for each \(i=1,\dotsc,n\):
\[
  k_{i}+\dotsb+k_{n}
  \leq
  0.
\]
Taking account that the coefficient \(\mathsf{C}_{\mathbi{k}}\) is zero unless:
\[
  k_{1}+\dotsb+k_{n}
  =
  0.
\]
It can be reformulated in saying that the coefficient \(\mathsf{C}_{\mathbi{k}}\) is zero unless for each \(i=1,\dotsc,n\):
\[
  k_{1}+\dotsb+k_{i}
  \geq
  0.
\]
This suggests to make a change of variables in order to deal with formal series.

In this work, we will not do this change of variables, but it is useful to introduce the length:
\[
  \ell(k_{1},k_{2},\dotsc,k_{n})
  =
  (k_{1})
  +
  (k_{1}+k_{2})
  \dotsb
  +
  (k_{1}+k_{2}+\dotsb+k_{n})
  =
  \sum_{i=1}^{n}
  (n-i)\,k_{i},
\]
that would correspond to the sum of the obtained exponents in the new variables.

We will also use the notation:
\[
  \ell(t_{1}^{k_{1}}\dotsm t_{n}^{k_{n}})
  \bydef
  \ell(k_{1},\dotsc,k_{n}),
\]
for the weighted degree of a monomial, and for \(l\in\N\), we will also write \(O(l)\) for a series that involves only monomials with weighted degree at least \(l\).
Using this notation, observe that:
\[
  \ell(t_i/t_j)=j-i.
\]

This basic observation will allow us to easily expand up to terms of degree \(3\) the factors appearing in the product:
\[
  \mathsf{C}(\underline{t})
  =
  \prod_{1\leq i<j\leq n}
  \frac{t_j-t_i}{t_j-2t_i}\;
  \prod_{2\leq i<j\leq n}
  \frac{t_j-2\,t_i}{t_j-2t_i+t_{i-1}}.
\]
\begin{Statement}
  For every integers \(i,j\) such that \(1\leq i<j\leq n\):
  \[
    \frac{t_j-t_i}{t_j-2\,t_i}
    =
    1
    +
    \frac{t_i}{t_j}
    +
    2
    \left(\frac{t_i}{t_j}\right)^2
    +
    O(3).
  \]
  and for every integers \(i,j\) such that \(2\leq i<j\leq n\):
  \[
    \frac{t_j-2\,t_i}{t_j-2\,t_i+t_{i-1}}
    =
    1
    -
    \frac{t_{i-1}}{t_j}
    +
    O(3).
  \]
\end{Statement}
\begin{proof}
  For every two integers \(i,j\) such that \(1\leq i<j\leq n\), the series expansion:
  \[
    \frac{t_j-t_i}{t_j-2\,t_i}
    =
    1
    +
    \frac{t_i}{t_j}
    \frac{1}{1-2\,\frac{t_i}{t_j}}
    =
    1
    +
    \frac{t_i}{t_j}
    \sum_{k\geq0}
    \left(2\,\frac{t_i}{t_j}\right)^k
  \]
  yields:
  \[
    \frac{t_j-t_i}{t_j-2\,t_i}
    =
    1
    +
    \frac{t_i}{t_j}
    +
    2
    \left(\frac{t_i}{t_j}\right)^2
    +
    O\bigl(3\,(j-i)\bigr).
  \]

  Next, for every two integers  \(i,j\) such that \(2\leq i<j\leq n\), the series expansion:
  \[
    \frac{t_j-2\,t_i}{t_j-2\,t_i+t_{i-1}}
    =
    1
    +
    \frac{-t_{i-1}}{t_j}
    \frac{1}{1-\frac{2\,t_i-t_{i-1}}{t_j}}
    =
    1
    -
    \frac{t_{i-1}}{t_j}
    \sum_{k\geq0}
    \left(\frac{2\,t_i-t_{i-1}}{t_j}\right)^k
  \]
  yields:
  \[
    \frac{t_j-2\,t_i}{t_j-2\,t_i+t_{i-1}}
    =
    1
    -
    \frac{t_{i-1}}{t_j}
    +
    O\bigl((j-i+1)+(j-i)\bigr).
  \]
\end{proof}
\begin{Statement}
  The expansion of \(\mathsf{C}\) with respect to the order:
  \[t_{1}\ll t_{2}\ll \dotsb\ll t_{n}\ll1\] 
  up to terms of weighted length \(\ell\) at least \(3\) is the following sum of terms having positive coefficients:
  \[
    \mathsf{C}(t_{1},\dotsc,t_{n})
    =
    1+
    \sum_{i=1}^{n-1}
    \frac{t_i}{t_{i+1}}
    +
    2\,
    \sum_{i=1}^{n-1}
    \frac{t_i^2}{t_{i+1}^2}
    +
    \sum_{i=1}^{n-2}
    \sum_{j=i+1}^{n-1}
    \frac{t_i}{t_{i+1}}
    \frac{t_j}{t_{j+1}}
    +
    O(3).
  \]
\end{Statement}
\begin{proof}
  Recall that \(\mathsf{C}(t_{1},\dotsc,t_{n})\) is by definition (the expansion of) the product:
  \[
    \mathsf{C}(t_{1},\dotsc,t_{n})
    =
    \prod_{1\leq i<j\leq n}
    \!\!\left(
    \frac{t_j-t_i}{t_j-2t_i}
    \right)\;
    \prod_{2\leq i<j\leq n}
    \!\!\left(
    \frac{t_j-2\,t_i}{t_j-2t_i+t_{i-1}}
    \right).
  \]

  By the lemma just above, the second product has the partial expansion:
  \[
    \prod_{2\leq i<j\leq n}
    \!\!\left(
    \frac{t_j-2\,t_i}{t_j-2t_i+t_{i-1}}
    \right)
    =
    \prod_{1\leq i\leq n-2}
    \left(
    1
    -
    \frac{t_{i}}{t_{i+2}}
    \right)
    +
    O\bigl(3\bigr).
  \]
  because for \(j\geq i+2\) the length \(j-i+1\geq3\).
  Thus:
  \[
    \prod_{2\leq i<j\leq n}
    \!\!\left(
    \frac{t_j-2\,t_i}{t_j-2t_i+t_{i-1}}
    \right)
    =
    1
    -
    \sum_{i=1}^{n-2}
    \frac{t_{i}}{t_{i+2}}
    +
    O\bigl(3\bigr).
  \]

  One has also:
  \[
    \prod_{1\leq i<j\leq n}
    \!\!\left(
    \frac{t_j-t_i}{t_j-2t_i}
    \right)
    =
    \prod_{1\leq i<j\leq n}
    \left(
    1
    +
    \frac{t_i}{t_j}
    +
    2
    \left(\frac{t_i}{t_j}\right)^2
    \right)
    +
    O(3).
    \\
    =
    \prod_{1\leq i\leq n-1}
    \left(
    1
    +
    \frac{t_i}{t_{i+1}}
    +
    2
    \left(\frac{t_i}{t_{i+1}}\right)^2
    \right)
    \prod_{1\leq i\leq n-2}
    \left(
    1
    +
    \frac{t_i}{t_{i+2}}
    \right)
    +
    O(3).
    \\
    =
    1+
    \sum_{i=1}^{n-1}
    \frac{t_i}{t_{i+1}}
    +
    2\,
    \sum_{i=1}^{n-1}
    \frac{t_i^2}{t_{i+1}^2}
    +
    \sum_{i=1}^{n-2}
    \sum_{j=i+1}^{n-1}
    \frac{t_i}{t_{i+1}}
    \frac{t_j}{t_{j+1}}
    +
    \sum_{i=1}^{n-2}
    \frac{t_{i}}{t_{i+2}}
    +
    O(3).
  \]

  We will consider separately \(j=2\), \(j=3\) and \(j\geq4\) and show that in all cases:
  \[
    \frac{t_j-t_{1}}{t_j-2t_{1}}
    \prod_{i=2}^{j-1}
    \frac{t_j-t_i}{t_j-2t_i+t_{i-1}}
    =
    1
    +
    \frac{t_{j-1}}{t_j}+2\left(\frac{t_{j-1}}{t_j}\right)^2
    +O(3).
  \]

  For the first case \(j=2\), there is only one term, and by the lemma just above:
  \[
    \frac{t_{2}-t_{1}}{t_{2}-2t_{1}}
    =
    1+\frac{t_{1}}{t_{2}}+2\left(\frac{t_{1}}{t_{2}}\right)^2+O(3).
  \]
  Hence, we are done with this case.

  For the second case \(j=3\), note that \(\ell\bigl((t_{1}/t_3)^2\bigr)=4\), thus:
  \[
    \frac{t_3-t_{1}}{t_3-2t_{1}}
    =
    1+\frac{t_{1}}{t_3}+O(3).
  \]
  There is only one supplementary term in the product, namely:
  \[
    \frac{t_3-t_{2}}{t_3-2\,t_{2}+t_{1}}
    =
    1
    +
    \frac{t_{2}}{t_3}
    -
    \frac{t_{1}}{t_3}
    +
    2\left(\frac{t_{2}}{t_3}\right)^2
    +
    O(3).
  \]
  By doing the product of these two terms, the coefficient of \(\frac{t_{1}}{t_3}\) cancels and all cross products of fractions are of weighted degree at least \(3\), hence:
  \[
    \frac{t_3-t_{1}}{t_3-2t_{1}}
    \frac{t_3-t_{2}}{t_3-2\,t_{2}+t_{1}}
    =
    1
    +
    \frac{t_{2}}{t_3}
    +
    2\left(\frac{t_{2}}{t_3}\right)^2
    +
    O(3).
  \]

  For the remaining case \(j\geq 4\), note that \(\ell(t_j/t_{1})\geq3\) thus:
  \[
    \frac{t_j-t_{1}}{t_j-2t_{1}}
    =
    1
    +
    O(3).
  \]
  In the same spirit, if \(j-i\geq3\):
  \[
    \frac{t_j-t_i}{t_j-2t_i+t_{i-1}}
    =
    1
    +
    O(3).
  \]
  The remaining factors of the product are obtained for \(i=j-1\):
  \[
    \frac{t_j-t_{j-1}}{t_j-2t_{j-1}+t_{j-2}}
    =
    1
    +
    \frac{t_{j-1}}{t_j}
    -
    \frac{t_{j-2}}{t_j}
    +
    2\left(\frac{t_{j-1}}{t_j}\right)^2
    +
    O(3),
  \]
  and for \(i=j-2\):
  \[
    \frac{t_j-t_{j-2}}{t_j-2t_{j-2}+t_{j-3}}
    =
    1
    +
    \frac{t_{j-2}}{t_j}
    +O(3).
  \]
  In the exact same way as above for \(j=3\), we obtain:
  \begin{multline*}
    \frac{t_j-t_{1}}{t_j-2t_{1}}\;
    \frac{t_j-t_{j-1}}{t_j-2t_{j-1}+t_{j-2}}\;
    \frac{t_j-t_{j-2}}{t_j-2t_{j-2}+t_{j-3}}\;
    \prod_{i=2}^{j-3}
    \frac{t_j-t_i}{t_j-2t_i+t_{i-1}}
    =\\
    1
    +
    \frac{t_{j-1}}{t_j}+2\left(\frac{t_{j-1}}{t_j}\right)^2
    +O(3).
  \end{multline*}

  It remains to state that the product of the obtained expressions for \(j=2,3,\dotsc,n\) is:
  \begin{multline*}
    \prod_{j=2}^{n}
    \biggl(
    1
    +
    \textcolor{gray}{\underbrace{\textcolor{black}{\left(\frac{t_{j-1}}{t_j}\right)}}_{\ell=1}}
    +
    2
    \textcolor{gray}{\underbrace{\textcolor{black}{\left(\frac{t_{j-1}}{t_j}\right)^2}}_{\ell=2}}
    +
    O(3)
    \biggr)
    =\\
    1+
    \left(
    \sum_{2\leq j\leq n}
    \biggl(
    \frac{t_{j-1}}{t_j}
    +
    2\,
    \frac{t_{j-1}^2}{t_j^2}
    \biggr)
    \right)
    +
    \left(
    \sum_{2\leq j_{1}<j_{2}\leq n}
    \frac{t_{j_{1}-1}}{t_{j_{1}}}
    \frac{t_{j_{2}-1}}{t_{j_{2}}}
    \right)
    +
    O(3).
  \end{multline*}
  This is because in the second parenthesis, that contains the cross products of fractions, the only terms that have degree least than \(3\) are the product of fractions of degree \(1\).
  By shifting all indices by \(-1\) one obtain the truncated expansion:
  \[
    \mathsf{C}(t_{1},\dotsc,t_{n})
    =
    1+
    \sum_{i=1}^{n-1}
    \frac{t_i}{t_{i+1}}
    +
    2\,
    \sum_{i=1}^{n-1}
    \frac{t_i^2}{t_{i+1}^2}
    +
    \sum_{i=1}^{n-2}
    \sum_{j=i+1}^{n-1}
    \frac{t_i}{t_{i+1}}
    \frac{t_j}{t_{j+1}}
    +
    O(3),
  \]
  as announced.
\end{proof}
Notice that in the new variables, this iterated Laurent series expansion of \(\mathsf{C}\) coincide with the usual multivariate Taylor expansion of th expression obtained from \(\mathsf{C}\).

We state that all terms of order less than \(3\) are non negative in such way that this ``Taylor'' series expansion  of \(\mathsf{C}\) allows to give a lower bound for the sum of positive contributions \(I_{0}^{+}\):
\begin{Statement}
  Under hypotheses \eqref{hyp:1} and \eqref{hyp:2}:
  \[
    I_{0}^{+}
    \geq
    2 \widetilde{I}_{0}.
  \]
\end{Statement}
\begin{proof}
  Recall that by definition, \(I_{0}^{+}\) is the following sum of positive contributions to the leading coefficient of \(I\):
  \[
    I_{0}^{+}
    =
    \sum_{\substack{
      k_{1}+\dotsb+k_{n}=0\\
      [\underline{t}^{\mathbi{k}}\,]\mathsf{C}(\underline{t})>0
    }}
    \bigl[t_{1}^{n-k_{1}}\dotsm t_{n}^{n-k_{n}}\bigr]
    \Bigl(
    f_{0}\bigl(\underline{t}\bigr)
    \Bigr)
    \bigl[t_{1}^{k_{1}}\dotsm t_{n}^{k_{n}}\bigr]
    \Bigl(
    \mathsf{C}\bigl(\underline{t}\bigr)
    \Bigr).
  \]
  Thanks to the multinomial formual, it is easy to compute that the appearing coefficients of \(f_{0}\) are:
  \[
    \bigl[
      t_{1}^{n-k_{1}}\dotsm t_{n}^{n-k_{n}}
    \bigr]
    \Bigl(
    f_{0}(t_{1},\dotsc,t_{n})\;
    \Bigr)
    =
    \frac{(a_{1}\dotsm a_{n})^{n}}{a_{1}^{k_{1}}\dotsm a_{n}^{k_{n}}}
    \frac{(n^{2})!}{(n-k_{1})!\dotsm(n-k_{n})!}.
  \]
  The coefficient \(\widetilde{I}_{0}\) is the same coefficient for \(\mathbi{k}=(0,\dotsc,0)\):
  \[
    \bigl[
      t_{1}^{n}\dotsm t_{n}^{n}
    \bigr]
    \Bigl(
    f_{0}(t_{1},\dotsc,t_{n})\;
    \Bigr)
    =
    (a_{1}\dotsm a_{n})^{n}
    \frac{(n^{2})!}{(n!)^{n}}.
  \]
  Thus:
  \[
    \bigl[
      t_{1}^{n-k_{1}}\dotsm t_{n}^{n-k_{n}}
    \bigr]
    \Bigl(
    f_{0}(t_{1},\dotsc,t_{n})\;
    \Bigr)
    =
    \frac{\widetilde{I}_{0}}{a_{1}^{k_{1}}\dotsm a_{n}^{k_{n}}}
    \frac{(n)!}{(n-k_{1})!}\dotsm\frac{(n)!}{(n-k_{n})!}.
  \]

  On the other hand, the above ``Taylor'' expansion of \(\mathsf{C}\) provide us with some points in the set
  \(\left\{k_{1}+\dotsb+k_{n}=0\colon[\underline{t}^{\mathbi{k}}\,]\mathsf{C}(\underline{t})>0\right\}\).
  Here is the list of the corresponding coefficients:
  \[
    \begin{array}{c|c|c|c|c|c}
      t^\mathbi{i}&1&t_i/t_{i+1}&t_i^2/t_{i+1}^2&t_i/t_{i+2}&t_it_j/t_{i+1}t_{j+1}\\\hline
      \mathrm{coeff}/\widetilde{I}_{0}
      &
      1
      &
      \frac{n\,a_{i+1}}{(n+1)\,a_i}
      &
      \frac{n(n-1)\,a_{i+1}^2}{(n+1)(n+2)\,a_i^2}
      &
      \frac{n\,a_{i+2}}{(n+1)\,a_i}
      &
      \frac{n^2\,a_{i+1}a_{j+1}}{(n+1)^2\,a_ia_j}
    \end{array}
  \]

  In conclusion:
  \begin{small}
    \begin{multline*}
      I_{0}^{+}
      \geq
      \widetilde{I}_{0}
      \Biggl(
      1+
      \sum_{i=1}^{n-1}
      \frac{n\,a_{i+1}}{(n+1)\,a_i}
      +
      2\,
      \sum_{i=1}^{n-1}
      \frac{n(n-1)\,a_{i+1}^2}{(n+1)(n+2)\,a_i^2}
      +
      \sum_{i=1}^{n-2}
      \frac{n\,a_{i+2}}{(n+1)\,a_i}
      +
      \sum_{i=1}^{n-3}
      \sum_{j=i+2}^{n-1}
      \frac{n^2\,a_{i+1}a_{j+1}}{(n+1)^2\,a_ia_j}
      \Biggr).
    \end{multline*}
  \end{small}
  For \(\underline{a}=\bigl(n^{n},\dotsc,n,1\bigr)\), one gets the sum:
  \[
    \frac{I_{0}^{+}}{\widetilde{I}_{0}}
    \geq
    \Biggl(
    1+
    \frac{(n-1)}{(n+1)}
    \frac{n}{n}
    +
    \frac{(n-1)^2}{(n+1)(n+2)}
    \frac{2\,n}{n^2}
    +
    \frac{(n-2)}{(n+1)}
    \frac{n}{n^2}
    +
    \frac{(n-3)(n-2)}{(n+1)^2}\frac{n^2}{2\,n^2}
    \Biggr),
  \]
  and we claim that, for \(n\geq5\), this quantity is more than \(2\).
\end{proof}

\subsection{Evaluation of \ensuremath{\lvert\mathsf{C}\rvert} and negative contributions}
\label{apx:C_bound}
Next, we control the negative contributions that could appear by multiplying \(\mathsf{A}\) by \(\mathsf{C}\):
\[
  I_{0}^{-}
  =
  -
  \sum_{\substack{
    k_{1}+\dotsb+k_{n}=0\\
    [\underline{t}^{\mathbi{k}}\,]\mathsf{C}(\underline{t})<0
  }}
  \bigl[t_{1}^{n-k_{1}}\dotsm t_{n}^{n-k_{n}}\bigr]
  \Bigl(
  f_{0}\bigl(\underline{t}\bigr)
  \Bigr)
  \bigl[t_{1}^{k_{1}}\dotsm t_{n}^{k_{n}}\bigr]
  \Bigl(
  \mathsf{C}\bigl(\underline{t}\bigr)
  \Bigr).
\]
We will use the transparent notation:
\[
  \mathsf{C}_{k_{1},\dotsc,k_{n}}
  \bydef
  \bigl[t_{1}^{k_{1}}\dotsm t_{n}^{k_{n}}\bigr]
  \mathsf{C}(t_{1},\dotsc,t_{n}).
\]
By triangular inequality, taking account of lemma \eqref{lem:fp} (for \(p=q=0\)), we get:
\[
  I_{0}^{-}
  \leq
  \sum_{\substack{
    k_{1}+\dotsb+k_{n}=0\\
    \mathsf{C}_{k_{1},\dotsc,k_{n}}<0
  }}
  \left(
  \frac{a_{1}}{2n\mu(\underline{a})}
  \right)^{-k_{1}}
  \dotsm
  \left(
  \frac{a_{n}}{2n\mu(\underline{a})}
  \right)^{-k_{n}}\;
  \widetilde{I}_{0}\;
  \abs{\mathsf{C}_{k_{1},\dotsc,k_{n}}}.
\]
Now, \(k_{1}+\dotsb+k_{n}=0\) thus:
\begin{equation}
  \label{eq:negativeContributions}
  I_{0}^{-}
  \leq
  \widetilde{I}_{0}
  \sum_{\substack{
    k_{1}+\dotsb+k_{n}=0\\
    \mathsf{C}_{k_{1},\dotsc,k_{n}}<0
  }}
  \abs{\mathsf{C}_{k_{1},\dotsc,k_{n}}}
  \left(
  \frac{1}{a_{1}}
  \right)^{k_{1}}
  \dotsm
  \left(
  \frac{1}{a_{n}}
  \right)^{k_{n}}.
\end{equation}
Recall that \(\abs{\mathsf{C}}(\underline{t})\) denotes the (convergent) iterated Laurent series generated by the absolute value of the coefficients of the series \(\mathsf{C}\):
\[
  \abs{\mathsf{C}}(t_{1},\dotsc,t_{n})
  \bydef
  \sum_{\mathbi{i}\in\Z^{n}}
  \abs{\mathsf{C}_{i_{1},\dotsc,i_{n}}}\;
  t_{1}^{i_{1}}\dotsm t_{n}^{i_{n}}.
\]
We will use the basic result:
\[
  \sum_{\substack{
    k_{1}+\dotsb+k_{n}=0\\
    \mathsf{C}_{k_{1},\dotsc,k_{n}}<0
  }}
  \abs{\mathsf{C}_{k_{1},\dotsc,k_{n}}}
  \left(
  \frac{1}{a_{1}}
  \right)^{k_{1}}
  \dotsm
  \left(
  \frac{1}{a_{n}}
  \right)^{k_{n}}
  =
  \frac{1}{2}
  \left(
  \abs{\mathsf{C}}\left(\frac{1}{a_{1}},\dotsc,\frac{1}{a_{n}}\right)
  -
  \mathsf{C}\left(\frac{1}{a_{1}},\dotsc,\frac{1}{a_{n}}\right)
  \right).
\]
\medskip

Now, we fix \(a_{1},\dotsc,a_{n}\), in order to establish the estimates:
\begin{Statement}
  \label{lem:C}
  Let \((a_i)_{i=1,\dotsc,n}\) be the decreasing geometric sequence 
  \[
    a_i
    \bydef
    (n)^{n-i}
    \qquad{\scriptstyle(i=1,\dotsc,n)},
  \]
  then for \(n\geq6\), the following inequalities hold:
  \[
    \frac{2}{3}\,\abs{\mathsf{C}}\biggl(\frac{1}{a_{1}},\dotsc,\frac{1}{a_{n}}\biggr)
    \leq
    \mathsf{C}\biggl(\frac{1}{a_{1}},\dotsc,\frac{1}{a_{n}}\biggr)
    \leq
    \abs{\mathsf{C}}\biggl(\frac{1}{a_{1}},\dotsc,\frac{1}{a_{n}}\biggr)
    \leq
    5.
  \]
\end{Statement}

Coming back to \eqref{eq:negativeContributions}, we get at once:
\[
  I_{0}^{-}
  \leq
  \widetilde{I}_{0}\,
  \frac{1}{2}
  \left(1-\frac{2}{3}\right)
  \abs{\mathsf{C}}\left(\frac{1}{a_{1}},\dotsc,\frac{1}{a_{n}}\right)\;
  \widetilde{I}_{0}
  \leq
  \frac{5}{6}\;
  \widetilde{I}_{0}.
\]

\begin{proof}[Proof of the inequalities]
  We insist on the fact that the coefficient of \(\mathsf{C}\) are very complicated. 
  In order to bypass this difficulty, we will use a \emph{majorant series} \(\widehat{\mathsf{C}}\) for \(\mathsf{C}\), in the sense that the Taylor series expansion of \(\widehat{\mathsf{C}}\) has only non negative coefficients, that are furthermore upper bounds for the absolute value of the corresponding Taylor coefficients of \(\mathsf{C}\):
  \[
    \abs{\mathsf{C}_{\mathbi{i}}}
    \leq
    \widehat{\mathsf{C}}_{\mathbi{i}}
    \qquad
    {\textstyle(\mathbi{i}\in\N^{n})}.
  \]
  Moreover, we will work in the domain of convergence of the series, in order to use their rational expressions.

  For each factor of 
  \[
    \mathsf{C}(t_{1},\dotsc,t_{n})
    =
    \prod_{1\leq i<j\leq n}
    \frac{t_{j}-t_{i}}{t_{j}-2t_{i}}\;
    \prod_{2\leq i<j\leq n}
    \frac{t_{j}-2\,t_{i}}{t_{j}-2t_{i}+t_{i-1}},
  \]
  it is easy to obtain such a majorant series, by replacing the series coefficients by their absolute values. 
  The first kind of fractions has only positive coefficients. Indeed:
  \[
    \frac{t_{j}-t_{i}}{t_{j}-2t_{i}}
    =
    1
    +
    \sum_{k\geq0}
    2^{k}
    \left(\frac{t_{i}}{t_{j}}\right)^{k+1}.
  \]
  The second kind of fractions has coefficients which sign is determined by the exponent of \(t_{i-1}\), indeed:
  \[
    \frac{t_{j}-2\,t_{i}}{t_{j}-2t_{i}+t_{i-1}}
    =
    1+\frac{(-t_{i-1})}{t_{j}}
    \sum_{k\geq l\geq0}
    (-1)^{l}\;
    2^{k-l}\binom{k}{l}\;
    t_{i-1}^{l}
    t_{i}^{k-l}
    t_{j}^{k}
  \]
  whereas, by changing the sign of \(t_{i-1}\):
  \[
    \frac{t_{j}-2\,t_{i}}{t_{j}-2t_{i}-t_{i-1}}
    =
    1+\frac{t_{i-1}}{t_{j}}
    \sum_{k\geq l\geq0}
    2^{k-l}\binom{k}{l}\;
    t_{i-1}^{l}
    t_{i}^{k-l}
    t_{j}^{k}.
  \]
  The majorant series is then obtained by taking the product of these pieces, that is:
  \[
    \widehat{\mathsf{C}}(t_{1},\dotsc,t_{n})
    =
    \prod_{1\leq i<j\leq n}
    \frac{t_{j}-t_{i}}{t_{j}-2t_{i}}\;
    \prod_{2\leq i<j\leq n}
    \frac{t_{j}-2\,t_{i}}{t_{j}-2t_{i}-t_{i-1}},
  \]
  The detail of this fact is elementary and left to the reader.

  Then:
  \[
    \widehat{\mathsf{C}}\left(\frac{1}{a_{1}},\dotsc,\frac{1}{a_{n}}\right)
    &=
    \prod_{1\leq i<j\leq n}
    \frac{a_{i}/a_{j}-1}{a_{i}/a_{j}-2}\;
    \prod_{2\leq i<j\leq n}
    \frac{a_{i}/a_{j}-2}{a_{i}/a_{j}-(2+a_{i}/a_{i-1})},
    \\
    &=
    \prod_{1\leq i<j\leq n}
    \frac{n^{j-i}-1}{n^{j-i}-2}\;
    \prod_{2\leq i<j\leq n}
    \frac{n^{j-i}-2}{n^{j-i}-(2+1/n)},
    \\
    &=
    \prod_{k=1}^{n-1}
    \left(\frac{n^{k}-1}{n^{k}-2}\right)^{n-k}\;
    \prod_{k=1}^{n-1}
    \left(\frac{n^{k}-2}{n^{k}-(2+1/n)}\right)^{n-1-k}\;
    \\
    &=
    \prod_{k=1}^{n-1}
    \left(
    \frac{\left(n^k-1\right)^{n-k}}
    {\left(n^k-2\right)\left(n^k-(2+1/n)\right)^{n-1-k}}
    \right)
  \]
  and similarly;
  \[
    \mathsf{C}\left(\frac{1}{a_{1}},\dotsc,\frac{1}{a_{n}}\right)
    =
    \prod_{k=1}^{n-1}
    \left(
    \frac{\left(n^{k}-1\right)^{n-k}}
    {\left(n^{k}-2\right)\left(n^{k}-(2-1/n)\right)^{n-1-k}}
    \right).
  \]

  This two products have many common terms and their difference is:
  \[
    \frac
    {\widehat{\mathsf{C}}\Bigl(\frac{1}{a_{1}},\dotsc,\frac{1}{a_{n}}\Bigr)}
    {\mathsf{C}\Bigl(\frac{1}{a_{1}},\dotsc,\frac{1}{a_{n}}\Bigr)}
    =
    \prod_{k=1}^{n-1}
    \frac{\left(n^{k}-2+1/n\right)^{n-1-k}}
    {\left(n^{k}-2-1/n\right)^{n-1-k}}
    =
    \prod_{k=1}^{n-1}
    \left(
    1+
    \frac{2}
    {n^{k+1}-2n-1)}
    \right)^{n-(k+1)}
  \]
  That is, after a shift of \(k\) by \(1\):
  \[
    \frac
    {\widehat{\mathsf{C}}\Bigl(\frac{1}{a_{1}},\dotsc,\frac{1}{a_{n}}\Bigr)}
    {\mathsf{C}\Bigl(\frac{1}{a_{1}},\dotsc,\frac{1}{a_{n}}\Bigr)}
    =
    \prod_{k=2}^{n}
    \left(
    1+
    \frac
    {2}
    {n^{k}-2n-1}
    \right)^{n-k}.
  \]

  In order to prove that for \(n\geq6\), this quantity is less than \(\frac{3}{2}\), consider the two first cases \(n=6,7\) with a computer algebra system, and then use a rough estimation for large enough cases, \textit{e.g.}:
  \[
    \ln\frac{\widehat{\mathsf{C}}}{\mathsf{C}}
    \leq
    \sum_{k=2}^{n}
    \frac{2(n-k)}{n^{k}-2n-1}
    \leq
    2
    \sum_{k=2}^{n}
    \frac{(n-2)}{(n-2)^{k}}
    \leq
    \frac{2}
    {(n-3)},
  \]
  that yields the upper bound by exponentiation, for \(n\geq8\).

  We estimate \(\widehat{\mathsf{C}}\) in the same way. We consider the cases \(7\leq n\leq 11\) with a computer algebra system and then one has: 
  \[
    \widehat{\mathsf{C}}\left(\frac{1}{a_{1}},\dotsc,\frac{1}{a_{n}}\right)
    =
    \prod_{k=1}^{n-1}
    \left(
    \frac{\left(n^k-1\right)^{n-k}}
    {\left(n^k-2\right)\left(n^k-(2+1/n)\right)^{n-1-k}}
    \right)
    \leq
    \prod_{k=1}^{n-1}
    \left(
    \frac{n^k-1}
    {n^k-(2+1/n)}
    \right)^{n-k}
  \]
  That is:
  \[
    \widehat{\mathsf{C}}\left(\frac{1}{a_{1}},\dotsc,\frac{1}{a_{n}}\right)
    \leq
    \prod_{k=1}^{n-1}
    \left(
    1+
    \frac{n+1}
    {n^{k+1}-2n-1}
    \right)^{n-k}
    =
    \prod_{k=2}^{n}
    \left(
    1+
    \frac{n+1}
    {n^{k}-2n-1}
    \right)^{n+1-k}
  \]
  We infer that:
  \[
    \ln(\widehat{\mathsf{C}})
    \leq
    \sum_{k=2}^{n}
    \frac{(n+1-k)(n+1)}
    {n^{k}-2n-1}
    \leq
    n^2
    \sum_{k\geq2}
    \left(\frac{1}{n-2}\right)^k
    \leq
    \frac{n^2}{(n-2)(n-3)}.
  \]
  For \(n\geq12\), by exponentiation:
  \[
    \widehat{\mathsf{C}}
    \left(
    \frac{1}{a_{1}},
    \dotsc,
    \frac{1}{a_{n}}
    \right)
    \leq
    \exp\left(\frac{n^2}{(n-2)(n-3)}\right)
    \leq 
    5.
  \]
\end{proof}

\subsection{Estimation of the slope}
\label{apx:delta}
Lastly, we study the slope of \(I_{0}\) with respect to \(\delta\):
\[
  I_{0}'
  =
  \sum_{k_{1}+\dotsb+k_{n}=0}
  \bigl[
    t_{1}^{n-k_{1}}\dotsm t_{n}^{n-k_{n}}
  \bigr]
  \Bigl(
  f_{1}(t_{1},\dotsc,t_{n})\,
  \mathfrak{s}_{1}(t_{1},\dotsc,t_{n})
  \Bigr)\;
  \bigl[
    t_{1}^{k_{1}}\dotsm t_{n}^{k_{n}}
  \bigr]
  \Bigl(
  \mathsf{C}(t_{1},\dotsc,t_{n})
  \Bigr).
\]
By triangular inequality, taking account of lemma \eqref{lem:fp} (for \(p=1\) and \(q=0\)), we get:
\[
  I_{0}^{-}
  \leq
  \widetilde{I}_{1}\;
  \sum_{k_{1}+\dotsb+k_{n}=0}
  \left(
  \frac{a_{1}}{2n\mu(\underline{a})}
  \right)^{-k_{1}}
  \dotsm
  \left(
  \frac{a_{n}}{2n\mu(\underline{a})}
  \right)^{-k_{n}}\;
  \abs{\mathsf{C}_{k_{1},\dotsc,k_{n}}}
  =
  \widetilde{I}_{1}\ 
  \abs{\mathsf{C}}
  \left(
  \frac{1}{a_{1}},\dotsc,\frac{1}{a_{n}}
  \right).
\]

On the other hand, in \eqref{cor:approximatedFujiwaraBound}, we have established that \(\widetilde{\lambda}=2\,\widetilde{I}_{1}/\widetilde{I}_{0}\), hence:
\[
  \abs{I_{0}'}
  &\leq
  \frac{\widetilde{\lambda}(\underline{a})}{2}\;
  \abs{\mathsf{C}}
  \left(
  \frac{1}{a_{1}},\dotsc,\frac{1}{a_{n}}
  \right)\
  \widetilde{I}_{0}
  \leq
  \frac{5\widetilde{\lambda}(\underline{a})}{2}
  \widetilde{I}_{0}.
\]

\newpage
\section{Estimation of the non-leading coefficients}
\label{apx:otherCoeffs}
\begin{itshape}
  In this appendix, we estimate the other coefficients. We also show that the term \(\mathsf{B}\) is always negligible.
\end{itshape}

\subsection{Estimation of the coefficient of \ensuremath{d^{n-1}} in \ensuremath{I}}
\label{apx:dn-1}
Next, we consider the coefficient of \(d^{n-1}\) in \(I\), that is the Cauchy product coefficient:
\[
  -I_{1}
  =
  \bigl[
    h^{n}t_{1}^{n}\dotsm t_{n}^{n}
  \bigr]
  \Bigl(
  \mathsf{A}_{1}(t_{1},\dotsc,t_{n})\;
  \mathsf{B}(t_{1},\dotsc,t_{n})\;
  \mathsf{C}(t_{1},\dotsc,t_{n})
  \Bigr).
\]
The term \(\mathsf{C}(t_{1},\dotsc,t_{n})\) does not involve the variable \(h\), thus:
\[
  I_{1}
  =
  \bigl[t_{1}^{n}\dotsm t_{n}^{n}\bigr]
  \biggl(
  \bigl[h^{n}\bigr]
  \Bigl(
  \mathsf{A}_{1}(t_{1},\dotsc,t_{n})\,
  \mathsf{B}(t_{1},\dotsc,t_{n})
  \Bigr)
  \mathsf{C}(t_{1},\dotsc,t_{n})
  \biggr).
\]

By lemma \eqref{eq:coeffA} above, the first factor of this product:
\[
  \mathsf{A}_{1}(t_{1},\dotsc,t_{n})
  =
  \bigl[d^{n-1}\bigr]\mathsf{A}(t_{1},\dotsc,t_{n})
\] 
is the following polynomial multiple of \(h^{n-1}\):
\begin{multline*}
  \mathsf{A}_{1}(t_{1},\dotsc,t_{n})
  =
  h^{n}
  \Bigl(
  \delta\,(n+1)\,
  f_{1}(\underline{t})\;
  \mathfrak{s}_{1}(\underline{t})
  -
  \delta\,
  f_{2}(\underline{t})\;
  \mathfrak{s}_{2}(\underline{t})
  \Bigr)
  \\
  +
  h^{n-1}
  \Bigl(
  f_{0}(\underline{t})\;
  \mathfrak{s}_{1}(\underline{t})
  -
  \delta\,
  f_{1}(\underline{t})\;
  \mathfrak{s}_{2}(\underline{t})
  \Bigr),
\end{multline*}
with the notation of \eqref{eq:coeffA}. Here we take account that: 
\[
  \alpha_{0}=1,\;
  \alpha_{1}=\delta(n+1),\;
  \text{ and }\;
  \beta_{1}=\beta_{2}=\delta.
\]

Once multiplied by
\(
\mathsf{A}_{1}
(t_{1},\dotsc,t_{n})
=
O\bigl(h^{n-1}\bigr)
\)
any multiple of \(h^2\) appearing in the expansion of \(\mathsf{B}\) would increase too much the exponent of \(h\) in order to reach the aimed monomial \(h^{n}t_{1}^{n}\dotsm t_{n}^{n}\).
Hence, in our computation, one can replace the iterated Laurent series:
\[
  \mathsf{B}(t_{1},\dotsc,t_{n})
  =
  \sum_{j_{1},\dotsc,j_{n}\geq0}
  \binom{n+j_{1}}{n}
  \dotsm
  \binom{n+j_{n}}{n}
  \;
  \frac{(-h)^{j_{1}+\dotsb+j_{n}}}
  {t_{1}^{j_{1}}\dotsm t_{n}^{j_{n}}},
\]
by its truncation up to polynomial multiples of \(h^2\):
\[
  \mathsf{B}(t_{1},\dotsc,t_{n})
  =
  1-
  \sum_{i=1}^{n}(n+1)
  \frac{h}{t_i}+O(h^2).
\]
The coefficient of \(h^{n}\) becomes:
\begin{multline*}
  \bigl[h^{n}\bigr]
  \Bigl(
  \mathsf{A}_{1}(t_{1},\dotsc,t_{n})\,
  \mathsf{B}(t_{1},\dotsc,t_{n})
  \Bigr)
  =
  \\
  \delta\,(n+1)\,
  f_{1}(\underline{t})\;
  \mathfrak{s}_{1}(\underline{t})
  -
  \delta\,
  f_{2}(\underline{t})\;
  \mathfrak{s}_{2}(\underline{t})
  \\-
  (n+1)
  \Bigl(
  f_{0}(\underline{t})\;
  \mathfrak{s}_{1}(\underline{t})
  -
  \delta\,
  f_{1}(\underline{t})\;
  \mathfrak{s}_{2}(\underline{t})
  \Bigr)
  \sum_{i=1}^{n}\frac{1}{t_{i}}.
\end{multline*}

Recall the convention:
\[
  \mathsf{C}_{k_{1},\dotsc,k_{n}}
  =
  \bigl[
    t_{1}^{k_{1}}\dotsm t_{n}^{k_{n}}
  \bigr]
  \Bigl(
  \mathsf{C}(t_{1},\dotsc,t_{n})
  \Bigr).
\]
Then:
\[
  \bigl[
    t_{1}^{k_{1}}\dotsm t_{n}^{k_{n}}
  \bigr]
  \left(
  \frac{1}{t_{i}}\mathsf{C}(t_{1},\dotsc,t_{n})
  \right)
  =
  \mathsf{C}_{k_{1},\dotsc,k_{i-1},k_{i}+1,k_{i+1},\dotsc,k_{n}}.
\]
Hence:
\begin{small}
  \begin{multline*}
    -I_{1}
    =
    \delta
    \sum_{k_{1}+\dotsb+k_{n}=0}
    \bigl[
      t_{1}^{n-k_{1}}\dotsm t_{n}^{n-k_{n}}
    \bigr]
    \Bigl(
    (n+1)\,
    f_{1}(\underline{t})\;
    \mathfrak{s}_{1}(\underline{t})
    -
    f_{2}(\underline{t})\;
    \mathfrak{s}_{2}(\underline{t})
    \Bigr)\;
    \mathsf{C}_{k_{1},\dotsc,k_{n}}
    \\-
    (n+1)
    \sum_{k_{1}+\dotsb+k_{n}=-1}
    \bigl[
      t_{1}^{n-k_{1}}\dotsm t_{n}^{n-k_{n}}
    \bigr]
    \Bigl(
    f_{0}(\underline{t})\;
    \mathfrak{s}_{1}(\underline{t})
    -
    \delta\,
    f_{1}(\underline{t})\;
    \mathfrak{s}_{2}(\underline{t})
    \Bigr)\;
    \bigl(\mathsf{C}_{k_{1}+1,\dotsc,k_{n}}+\dotsb+\mathsf{C}_{k_{1},\dotsc,k_{n}+1}\bigr)
  \end{multline*}
\end{small}
By triangular inequality, applying four times lemma \eqref{lem:fp}:
\begin{small}
  \begin{multline*}
    \abs{I_{1}}
    \leq
    \delta\;\bigl((n+1)\;\widetilde{I}_{1}+\widetilde{I}_{2}\bigr)
    \sum_{k_{1}+\dotsb+k_{n}=0}
    \left(\frac{1}{a_{1}}\right)^{k_{1}}
    \dotsm
    \left(\frac{1}{a_{n}}\right)^{k_{n}}
    \;
    \abs{\mathsf{C}}_{k_{1},\dotsc,k_{n}}
    \\
    +
    (n+1)(\widetilde{I}_{1}+\delta\,\widetilde{I}_{2})
    \sum_{k_{1}+\dotsb+k_{n}=-1}
    \frac{1}{2n\mu(\underline{a})}
    \left(\frac{1}{a_{1}}\right)^{k_{1}}
    \dotsm
    \left(\frac{1}{a_{n}}\right)^{k_{n}}
    \;
    \bigl(\abs{\mathsf{C}}_{k_{1}+1,\dotsc,k_{n}}+\dotsb+\abs{\mathsf{C}}_{k_{1},\dotsc,k_{n}+1}\bigr).
  \end{multline*}
\end{small}
Now the sum in the second line expands as follows in \(n\) sums: 
\begin{multline*}
  \sum_{k_{1}+\dotsb+k_{n}=-1}
  \left(\frac{1}{a_{1}}\right)^{k_{1}}
  \dotsm
  \left(\frac{1}{a_{n}}\right)^{k_{n}}
  \;
  \bigl(\abs{\mathsf{C}}_{k_{1}+1,\dotsc,k_{n}}+\dotsb+\abs{\mathsf{C}}_{k_{1},\dotsc,k_{n}+1}\bigr)
  =\\
  \sum_{(k_{1}+1)+\dotsb+k_{n}=0}
  \left[
    \left(\frac{1}{a_{1}}\right)^{(k_{1}+1)-1}
    \dotsm
    \left(\frac{1}{a_{n}}\right)^{k_{n}} \;
    \abs{\mathsf{C}}_{k_{1}+1,\dotsc,k_{n}}
  \right]+
  \\
  \dotsb
  \\
  +
  \sum_{k_{1}+\dotsb+(k_{n}+1)=0}
  \left[
    \left(\frac{1}{a_{1}}\right)^{k_{1}}
    \dotsm
    \left(\frac{1}{a_{n}}\right)^{(k_{n}+1)-1} \;
    \abs{\mathsf{C}}_{k_{1},\dotsc,k_{n}+1}
  \right],
\end{multline*}
and it suffices to make the appropriate shift of indices in each of these sums in order to obtain:
\begin{multline*}
  \sum_{k_{1}+\dotsb+k_{n}=-1}
  \left(\frac{1}{a_{1}}\right)^{k_{1}}
  \dotsm
  \left(\frac{1}{a_{n}}\right)^{k_{n}}
  \;
  \bigl(\abs{\mathsf{C}}_{k_{1}+1,\dotsc,k_{n}}+\dotsb+\abs{\mathsf{C}}_{k_{1},\dotsc,k_{n}+1}\bigr)
  =\\
  a_{1}
  \abs{\mathsf{C}}
  \left(
  \frac{1}{a_{1}},
  \dotsc,
  \frac{1}{a_{n}}
  \right)
  +
  \dotsb
  +
  a_{n}
  \abs{\mathsf{C}}
  \left(
  \frac{1}{a_{1}},
  \dotsc,
  \frac{1}{a_{n}}
  \right).
\end{multline*}
The later sum is smaller than \(\mu(\underline{a})\abs{\mathsf{C}} \bigl( \frac{1}{a_{1}}, \dotsc, \frac{1}{a_{n}} \bigr)\) because by definition of \(\mu(\underline{a})\):
\[
  \mu(a_{1},\dotsc,a_{n})
  =
  1\,a_{1}
  +
  \dotsb
  +
  n\,a_{n}.
\]

One can now obtain, as desired, a constant multiple of \(\widetilde{I}_{1}\) bounding the coefficient \(\abs{I_{1}}\) from above:
\[
  \abs{I_{1}}
  &\leq
  \left(
  \delta\;\bigl((n+1)\;\widetilde{I}_{1}+\widetilde{I}_{2}\bigr)
  +
  \frac{n+1}{2n}
  (\widetilde{I}_{1}+\delta\,\widetilde{I}_{2})
  \right)\;
  \abs{\mathsf{C}}
  \left(
  \frac{1}{a_{1}},
  \dotsc,
  \frac{1}{a_{n}}
  \right)
  \\
  &\leq
  \widetilde{I}_{1}
  \biggl(
  \delta\;(n+1+\widetilde{\lambda})
  +
  \frac{n+1}{2n}
  (1+\delta\,\widetilde{\lambda})
  \biggr)\;
  \abs{\mathsf{C}}
  \left(
  \frac{1}{a_{1}},
  \dotsc,
  \frac{1}{a_{n}}
  \right)
  \\
  \abs{I_{1}}
  &\leq
  \widetilde{I}_{1}
  \biggl(
  \frac{1+3\,\delta\widetilde{\lambda}}{2}
  +
  \frac{1+\delta\widetilde{\lambda}}{2n}
  +
  \delta(n+1)
  \biggr)
  \abs{\mathsf{C}}
  \left(
  \frac{1}{a_{1}},
  \dotsc,
  \frac{1}{a_{n}}
  \right).
\]

Notice that this estimate are true without assumption \eqref{hyp:1} and \eqref{hyp:2} and also for \(\kappa\neq n\).
Now, taking account of \(5\delta\widetilde{\lambda}\leq1\) and \(n\geq6\), one has:
\[
  \left(
  \frac{1+3\,\delta\widetilde{\lambda}}{2}
  +
  \frac{1+\delta\widetilde{\lambda}}{2n}
  +
  \delta(n+1)
  \right)
  \leq
  \biggl(
  \frac{8}{10}
  +
  \frac{1}{10}
  +
  \frac{n+1}{5\widetilde{\lambda}}
  \biggr)
  \leq
  1.
\]
Because:
\[
  \widetilde{\lambda}
  \geq
  2\,(n+1),
\]
one has finally (without using hypothesis \eqref{hyp:1}):
\[
  \abs{I_{1}}
  \leq
  \abs{\mathsf{C}}
  \left(
  \frac{1}{a_{1}},
  \dotsc,
  \frac{1}{a_{n}}
  \right)\;
  \widetilde{I}_{1}.
\]

\subsection{Estimation of the other coefficients}
\label{apx:Ip}
We consider the coefficients of the remaining monomials \(d^{n-p}=d^{n-2},\dotsc,d^{1},d^{0}\).
As we have already seen twice above, the coefficient of \(d^{n-p}\) is a polynomial multiple of \(h^{n-p}\). Because we seek the coefficient of \(h^{n}\), this allowed us to simplify the computation by replacing the iterated Laurent series:
\[
  \mathsf{B}(t_{1},\dotsc,t_{\kappa})
  =
  \sum_{j_{1},\dotsc,j_{\kappa}\geq0}
  \binom{n+j_{1}}{n}
  \dotsm
  \binom{n+j_{\kappa}}{n}
  \;
  \frac{(-h)^{j_{1}+\dotsb+j_{\kappa}}}
  {t_{1}^{j_{1}}\dotsm t_{\kappa}^{j_{\kappa}}}
\]
by its truncations \(B\simeq1+O(h)\) for \(p=0\) and \(B\simeq1-\sum_{i=1}^\kappa(n+1)\frac{h}{t_i}+O(h^2)\) for \(p=1\).
As \(p\) goes from \(0\) to \(n\), the number of terms that we have to consider will increase dramatically. 
Consequently, the precise estimation of the coefficient:
\[
  -I_{p}
  =
  \bigl[d^{n-p}h^{n}t_{1}^{n}\dotsm t_{\kappa}^{n}\bigr]
  \Bigl(
  \mathsf{A}(t_{1},\dotsc,t_{\kappa})\;
  \mathsf{B}(t_{1},\dotsc,t_{\kappa})\;
  \mathsf{C}(t_{1},\dotsc,t_{\kappa})
  \Bigr)
\]
is technically more and more involved.

This consideration quickly justify that we will not use anymore the expanded expression:
\[
  \mathsf{B}(t_{1},\dotsc,t_{\kappa})
  =
  \sum_{j_{1},\dotsc,j_{\kappa}\geq0}
  \binom{n+j_{1}}{n}
  \dotsm
  \binom{n+j_{\kappa}}{n}
  \;
  \frac{(-h)^{j_{1}+\dotsb+j_{\kappa}}}
  {t_{1}^{j_{1}}\dotsm t_{\kappa}^{j_{\kappa}}}
\]
as in the preceding estimations of \(I_{0}\) and \(I_{1}\), but that we will use a simpler expression instead.

In order to use the triangular inequality, we will forget about the signs of the coefficients of \(\mathsf{B}\) and use the corresponding series with non negative coefficients:
\[
  \abs{\mathsf{B}}(t_{1},\dotsc,t_{\kappa})
  =
  \sum_{j_{1},\dotsc,j_{\kappa}\geq0}
  \binom{n+j_{1}}{n}
  \dotsm
  \binom{n+j_{\kappa}}{n}
  \;
  \frac{(\mathbf{+}h)^{j_{1}+\dotsb+j_{\kappa}}}
  {t_{1}^{j_{1}}\dotsm t_{\kappa}^{j_{\kappa}}}
\]
and moreover we will manage to use its rational expression:
\[
  \abs{\mathsf{B}}(t_{1},\dotsc,t_{\kappa})
  =
  \prod_{i=1}^{\kappa}
  \frac{1}{\bigl(1-\frac{h}{t_{1}}\bigr)^{n+1}},
\]
the estimation of which is more easy. In a sense, we will not need anymore to know anything about the coefficients of \(\mathsf{B}\). 

This level of precision will be sufficient in view of the application of Fujiwara's bound:
\[
  d
  \geq
  2\;
  \max_{1\leq p\leq n}
  \bigg\lvert
  \frac{I_{p}}{I_0}
  \bigg\rvert^{\frac{1}{p}},
\]
because we see that for \(p\geq2\), the obtained ratio will be decreased by an exponent \(\frac{1}{p}\leq\frac{1}{2}\).

\begin{Statement}
  For \(p=0,1,\dotsc,n\) and any choice of the parameters \(a_{1},\dotsc,a_{\kappa}\), the absolute value of the coefficient of \(d^{n-p}\) in the polynomial \(I\) is bounded from above by:
  \[
    \abs{I_{p}}
    \leq
    \bigl(\widetilde{I}_{p}+\delta\widetilde{I}_{p+1}\bigr)\;
    \abs{\mathsf{B}}\left(\frac{2n\mu h}{a_{1}},\dotsc,\frac{2n\mu h}{a_{\kappa}}\right)\,
    \abs{\mathsf{C}}\left(\frac{1}{a_{1}},\dotsc,\frac{1}{a_{\kappa}}\right).
  \]
\end{Statement}
\begin{proof}
  Recall that \(I(d)\) is the Cauchy product coefficient:
  \[
    I
    =
    \bigl[h^{n}t_{1}^{n}\dotsm t_{\kappa}^{n}\bigr]
    \Bigl(
    \mathsf{A}(t_{1},\dotsc,t_{\kappa})\,
    \mathsf{B}(t_{1},\dotsc,t_{\kappa})\,
    \mathsf{C}(t_{1},\dotsc,t_{\kappa})
    \Bigr).
  \]
  Recall that:
  \[
    \mathsf{B}(t_{1},\dotsc,t_{\kappa})
    =
    \sum_{\mathbi{j}\in\N^{\kappa}}
    \mathsf{B}_{\mathbi{j}}\,
    \Biggl(\frac{h}{t_{1}}\Biggr)^{\!\!j_{1}}
    \dotsm
    \Biggl(\frac{h}{t_{\kappa}}\Biggr)^{\!\!j_{\kappa}}.
  \]
  and:
  \[
    \mathsf{C}(t_{1},\dotsc,t_{\kappa})
    =
    \sum_{k_{1}+\dotsb+k_{\kappa}=0}
    \mathsf{C}_{\mathbi{k}}\,
    t_{1}^{k_{1}}
    \dotsm
    t_{\kappa}^{k_{\kappa}}.
  \]
  The coefficient od \(d^{n-p}\) is thus:
  \[
    -I_{p}
    =
    \sum_{\mathbi{i}-\mathbi{j}+\mathbi{k}=0}
    \bigl[h^{n-q}\underline{t}^{n+\mathbi{i}}\bigr]
    \mathsf{A}_{p}(h,\underline{t})\,
    \mathsf{B}_{\mathbi{j}}
    \mathsf{C}_{\mathbi{k}}
    \pour{q=i_{1}+\dotsb+i_{\kappa}=j_{1}+\dotsb+j_{\kappa}}.
  \]

  Next, we expand this expression by using \eqref{eq:coeffA}, and we get:
  \begin{multline*}
    -I_{p}
    =
    \sum_{\mathbi{i}-\mathbi{j}+\mathbi{k}=0}
    \alpha_{p-q}\;
    \bigl[\underline{t}^{n+\mathbi{i}}\bigr]
    \Bigl(
    f_{p-q}(\underline{t})\,\mathfrak{s}_{\kappa-n+p}(\underline{t})
    \Bigr)\,
    \mathsf{B}_{\mathbi{j}}
    \mathsf{C}_{\mathbi{k}}
    +\\
    \sum_{\mathbi{i}-\mathbi{j}+\mathbi{k}=0}
    \beta_{p-q+1}\;
    \bigl[\underline{t}^{n+\mathbi{i}}\bigr]
    \Bigl(
    f_{p-q+1}(\underline{t})\,\mathfrak{s}_{\kappa-n+p+1}(\underline{t})
    \Bigr)\,
    \mathsf{B}_{\mathbi{j}}
    \mathsf{C}_{\mathbi{k}},
  \end{multline*}
  where again:
  \[
    q=i_{1}+\dotsb+i_{\kappa}=j_{1}+\dotsb+j_{\kappa}.
  \]

  Recall that \(\abs{\alpha_{p-q}}\leq1\) and \(\abs{\beta_{p-q+1}}\leq\delta\). By triangular inequality and using \eqref{lem:fp}:
  \begin{multline*}
    \lvert I_{p}\rvert
    \leq
    \widetilde{I}_{p}
    \sum_{\substack{
      j_{1}+\dotsb+j_{\kappa}=q\\
      k_{1}+\dotsb+k_{\kappa}=0
    }}
    1\;
    \left(\frac{1}{2n\mu}\right)^{q}
    \frac{a_{1}^{j_{1}}\dotsm a_{\kappa}^{j_{\kappa}}}
    {a_{1}^{k_{1}}\dotsm a_{\kappa}^{k_{\kappa}}}
    \abs{\mathsf{B}_{\mathbi{j}}}
    \abs{\mathsf{C}_{\mathbi{k}}}
    +\\
    \widetilde{I}_{p+1}
    \sum_{\substack{
      j_{1}+\dotsb+j_{\kappa}=q\\
      k_{1}+\dotsb+k_{\kappa}=0
    }}
    \delta\;
    \left(\frac{1}{2n\mu}\right)^{q}
    \frac{a_{1}^{j_{1}}\dotsm a_{\kappa}^{j_{\kappa}}}
    {a_{1}^{k_{1}}\dotsm a_{\kappa}^{k_{\kappa}}}
    \abs{\mathsf{B}_{\mathbi{j}}}
    \abs{\mathsf{C}_{\mathbi{k}}},
  \end{multline*}
  We factorize \(\widetilde{I}_{p}+\delta\,\widetilde{I}_{p+1}\) and separate \(\mathbi{j}\) and \(\mathbi{k}\) and obtain:
  \[
    \big\lvert
    I_{p}
    \big\rvert
    \leq
    \bigl(\widetilde{I}_{p}+\delta\,\widetilde{I}_{p+1}\bigr)\;\;
    \sum_{j_{1}+\dotsb+j_{\kappa}=q}
    \left(\frac{1}{2n\mu}\right)^{q}
    a_{1}^{j_{1}}\dotsm a_{\kappa}^{j_{\kappa}}
    \abs{\mathsf{B}_{\mathbi{j}}}\;
    \sum_{k_{1}+\dotsb+k_{\kappa}=0}
    \frac{1}{a_{1}^{k_{1}}\dotsm a_{\kappa}^{k_{\kappa}}}
    \abs{\mathsf{C}_{\mathbi{k}}}.
  \]
  The announced result is obtained by identifying the two appearing sums:
  \[
    \abs{I_{p}}
    \leq
    \bigl(\widetilde{I}_{p}+\delta\widetilde{I}_{p+1}\bigr)\;
    \abs{\mathsf{B}}\left(\frac{2n\mu h}{a_{1}},\dotsc,\frac{2n\mu h}{a_{\kappa}}\right)\;
    \abs{\mathsf{C}}\left(\frac{1}{a_{1}},\dotsc,\frac{1}{a_{\kappa}}\right).
  \]
\end{proof}

We have already estimated \(\abs{\mathsf{C}}\left(\frac{1}{a_{1}},\dotsc,\frac{1}{a_{\kappa}}\right)\).
In order to apply Fujiwara's bound,it remains to estimate \(\abs{\mathsf{B}}\bigl(\frac{2n\mu h}{a_{1}},\dotsc,\frac{2n\mu h}{a_{\kappa}}\bigr)\).
We will now see that ---\,without using hypothese \eqref{hyp:1} and \eqref{hyp:2}\,--- this term is bounded from above by the term of a converging (thus bounded) sequence.
\begin{Statement}
  For any integer \(\kappa\in\N\), for any choice of the parameters \((a_{1},\dotsc,a_{\kappa})\in\N^\kappa\) and for 
  \(\mu=(1\,a_{1}+\dotsb+\kappa a_{\kappa})\).
  \[
    \abs{\mathsf{B}}
    \left(
    \frac{2n\mu h}{a_{1}},
    \dotsc,
    \frac{2n\mu h}{a_{\kappa}}
    \right)
    \leq
    \biggl(
    \frac{2\,n}{2\,n-1}
    \biggr)^{n+1}
    \quad
    {\scriptstyle
      \left(
      \to
      e^{1/2}
    \right)},
  \]
\end{Statement}
\begin{proof}
  Recall that:
  \[
    \abs{\mathsf{B}}(t_{1},\dotsc,t_{\kappa})
    =
    \sum_{j_{1},\dotsc,j_{\kappa}\geq0}
    \binom{n+j_{1}}{n}
    \dotsm
    \binom{n+j_{\kappa}}{n}\,
    \frac{h^{j_{1}+\dotsb+j_{\kappa}}}
    {t_{1}^{j_{1}}\dotsm t_{\kappa}^{j_{\kappa}}}.
  \]
  Thus:
  \[
    \abs{\mathsf{B}}
    \left(
    \frac{2n\mu h}{a_{1}},
    \dotsc,
    \frac{2n\mu h}{a_{\kappa}}
    \right)
    =
    \sum_{j_{1},\dotsc,j_{\kappa}\geq0}
    \binom{n+j_{1}}{n}
    \dotsm
    \binom{n+j_{\kappa}}{n}\,
    \left(\frac{a_{1}}{2n\mu}\right)^{j_{1}}
    \dotsm 
    \left(\frac{a_{\kappa}}{2n\mu}\right)^{j_{\kappa}}
  \]
  Now \(2n\mu>a_{i}\) for \(i=1,\dotsc,\kappa\), hence the left hand side expression is:
  \[
    \abs{\mathsf{B}}
    \left(
    \frac{2n\mu h}{a_{1}},
    \dotsc,
    \frac{2n\mu h}{a_{\kappa}}
    \right)
    =
    \prod_{i=1}^{\kappa}
    \left(1-\frac{a_i}{\mu}\frac{1}{2n}\right)^{-(n+1)},
  \]
  and by taking the logarithm:
  \[
    \ln
    \abs{\mathsf{B}}
    \left(
    \frac{2n\mu h}{a_{1}},
    \dotsc,
    \frac{2n\mu h}{a_{\kappa}}
    \right)
    =
    (n+1)
    \sum_{i=1}^{\kappa}
    \left\lvert
    \ln\left(1-\frac{a_i}{\mu}\frac{1}{2n}\right)
    \right\rvert.
  \]

  Next, use the concavity of the logarithm function and the fact that:
  \[
    0
    \leq
    a_i
    \leq
    \mu(\underline{a})
    =
    1\,a_{1}+\dotsb+\kappa\,a_{\kappa}
    \qquad
    {\scriptstyle (i=1,\dotsc,\kappa)}
  \]
  in order to obtain:
  \[
    \left\lvert
    \ln\left(1-\frac{a_{i}}{\mu}\,\frac{1}{2n}\right)
    \right\rvert
    \leq
    \frac{a_i}{\mu}
    \left\lvert
    \ln\left(1-\frac{1}{2\,n}\right)
    \right\rvert
    =
    \frac{a_i}{\mu}
    \ln\left(\frac{2\,n}{2\,n-1}\right)
    \qquad
    \pour{i=1,\dotsc,\kappa},
  \]
  and sum these inequalities in order to get the upper bound:
  \[
    \ln
    \abs{\mathsf{B}}
    \left(
    \frac{2n\mu h}{a_{1}},
    \dotsc,
    \frac{2n\mu h}{a_{\kappa}}
    \right)
    &\leq
    \frac{(a_{1}+\dotsb+a_{\kappa})}{\mu(\underline{a})}\,
    (n+1)
    \ln\left(\frac{2\,n}{2\,n-1}\right)
    \\&\leq
    (n+1)
    \ln\left(\frac{2\,n}{2\,n-1}\right).
  \]

  Finally, use the increasingness of the exponential function in order to prove the announced result:
  \[
    \abs{\mathsf{B}}
    \left(
    \frac{2n\mu h}{a_{1}},
    \dotsc,
    \frac{2n\mu h}{a_{\kappa}}
    \right)
    \leq
    \left(\frac{2\,n}{2\,n-1}\right)^{n+1}.
  \]
\end{proof}

\bibliographystyle{amsplain}
\bibliography{$LATEX/these}
\end{document}